\magnification 1195    
\overfullrule0pt

\def\sqr#1#2{{\vcenter{\hrule height.#2pt              
     \hbox{\vrule width.#2pt height#1pt\kern#1pt
     \vrule width.#2pt}
     \hrule height.#2pt}}}
\def\square{\mathchoice\sqr{5.5}4\sqr{5.0}4\sqr{4.8}3\sqr{4.8}3}
\def\qed{\hskip4pt plus1fill\ $\square$\par\medbreak}

\font\contrm=cmr9 at 8pt
\input epsf.sty
\def\d{{\rm deg}}

\def\cB{{\cal B}}

\def\cD{{\cal D}}

\def\cG{{\cal G}}
\def\cH{{\cal H}}

\def\cV{{\cal V}}
\def\cW{{\cal W}}

\def\C{{\bf C}}

\def\R{{\bf R}}

\centerline{\bf A Symbolic Characterization of the Horseshoe Locus}

\centerline{\bf in the H\'enon Family}
\bigskip

\centerline{Eric Bedford  and John Smillie}

\bigskip
{\obeylines \contrm
0.  Introduction
1.  Definition and properties of crossed maps
2.  Codings of Orbits
3.  The 3-Box System
4.  Real, Crossed mappings; The Disk Property
5.  The 3-Box System:  Real Case
6.  Expansion
7.  External Rays}

\bigskip\noindent{\bf 0. Introduction.}  The  H\'enon family has been studied as a model family both for the interesting dynamical behavior of particular maps  in the family (as in the work of Benedicks and Carleson and others on the existence of strange attractors) as well as for the way in which the dynamics varies with the parameter (as in the work of Newhouse on persistent instability).  The H\'enon family contains parameters that produce hyperbolic horseshoes and the boundary of this horseshoe region can be seen as representing one model for loss of hyperbolicity. Dynamical behavior of this family on the boundary of the horseshoe locus has been studied by several authors (see, for instance, [CLR1,2],  [H1,2],  [T]).

\medskip

\epsfysize=1.5in
\centerline{ \epsfbox{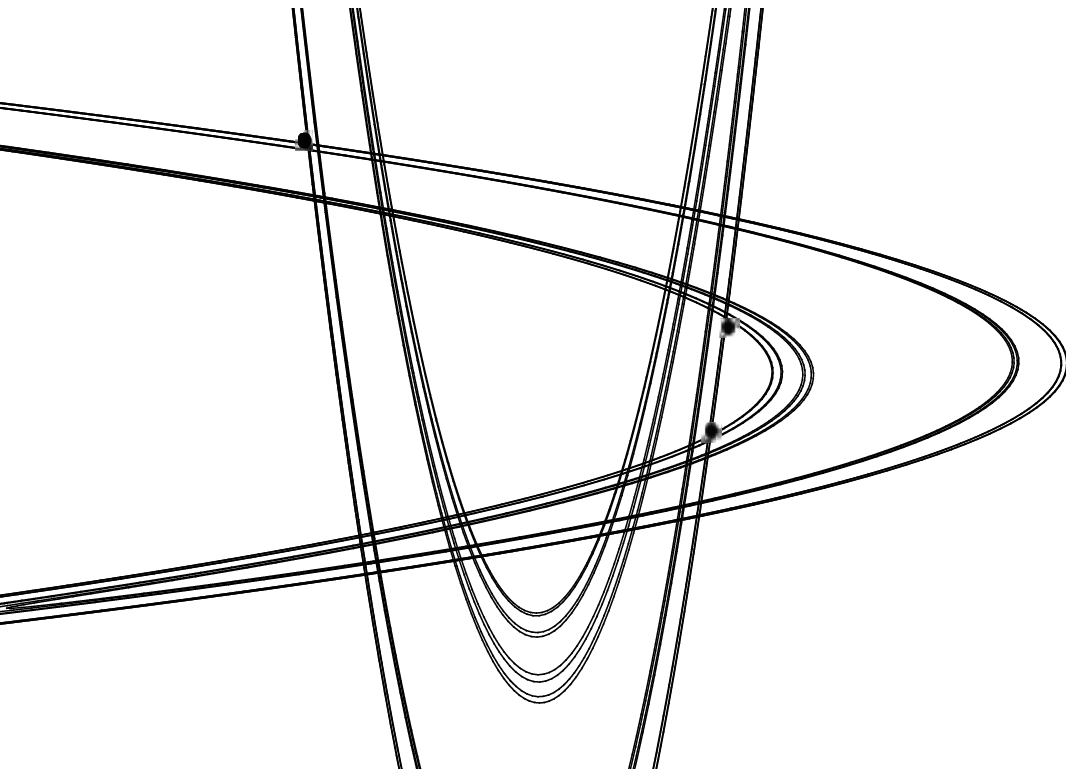}  }

\centerline{Figure 0.1.  Horseshoe for $f_{c,\delta}(x,y)= (c + \delta y - x^2, -x)$ with $c=6.0$, $\delta = 0.8$. }
\medskip

It will be worthwhile reviewing the one variable theory of the quadratic map since it gives a good starting point for describing our coding scheme.  If we fix a parameter and consider a single map then we can assign to each point in the dynamical interval an itinerary which encodes its position with respect to the critical point. The possible itineraries of all points are determined by the itinerary of the critical value. As the parameter changes we can consider the change of the itinerary of the critical value. This  proves to be a very useful dynamically significant function on parameter space. 
Kneading theory cannot be directly applied to the H\'enon family because there is no critical point, but the attempt to extend kneading theory to the H\'enon family leads to the theory of the pruning front initiated by Cvitanovic. This theory leads to significant challenges, both theoretical (see De Carvalho and Hall [DCH]) and computational  (see Hagiwara and Shudo [HS]).

\medskip
\epsfysize=1.4in
\centerline{ \epsfbox{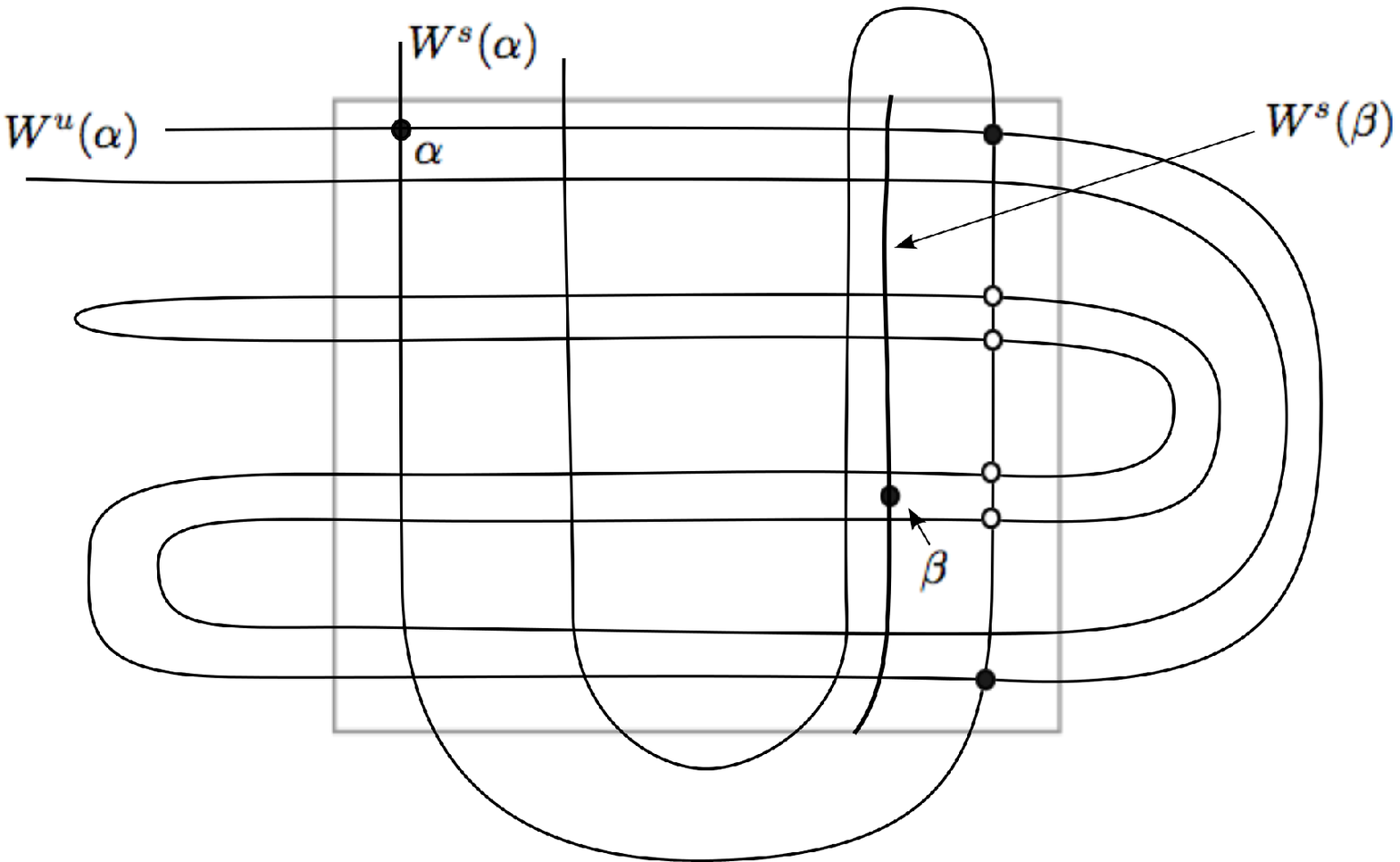}  }

\centerline{Figure 0.2.  Horseshoe: 2 points $\bullet$ are $\overline 0 1 \overline 0$; 4 points $\circ$ are $\overline0101\overline 0$.} 

\medskip

\medskip
\epsfysize=1.2in
\centerline{ \epsfbox{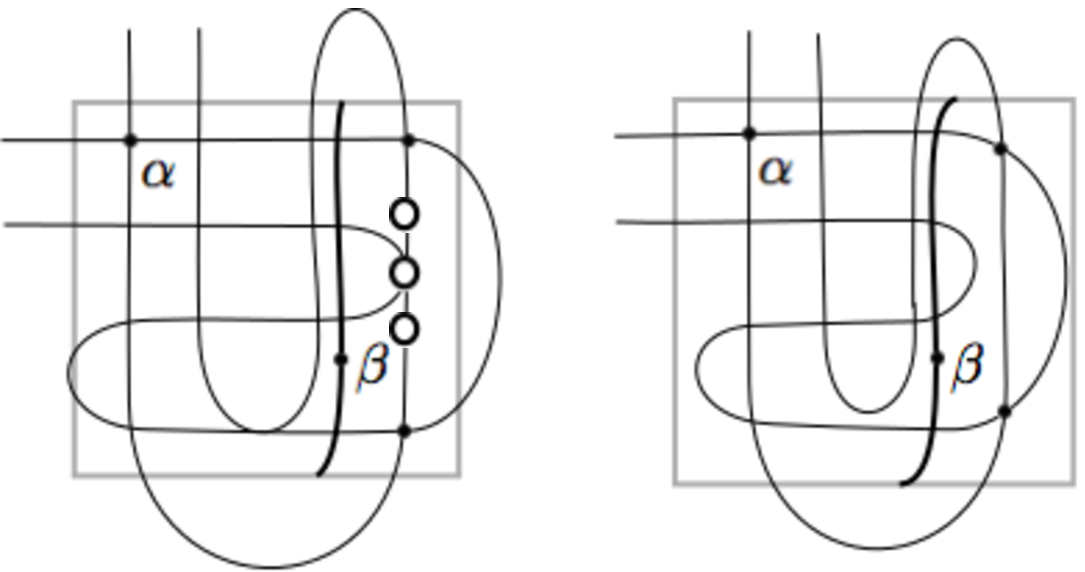}  }

\centerline{Figure 0.3.  Degeneration of the horseshoe:  2 points $\circ$ have come together. } 

\medskip

The iconic picture of the horseshoe map is the homoclinic tangle shown in Figure~0.1.  The marked point on the upper left is the unique one-sided saddle fixed point; the vertical/horizontal directions are the stable/unstable directions, respectively.  The unstable manifolds have a marked bend on the right, and there is an ordering of these loops.  The innermost loop is indicated by two marked points on the right.  The naive picture of degeneration is that the inner bend pulls to the left and, when the two marked points disappear, the horseshoe is destroyed.    It is this intuitive picture that we want to justify rigorously.  

Figure~0.2 is a stylized version of Figure~0.1.  There are two saddle fixed points, which are marked $\alpha$ and $\beta$, and an arc of the stable manifold $W^s(\beta)$ is indicated.  This stable arc cuts the box into two pieces, and we use it to define a ``right-left''  coding of an itinerary, with ``0'' to the left and ``1'' on the right, which serves as a rough analogue of the kneading sequence in one variable. Our coding has some advantages and some disadvantages when compared to the kneading sequence. Its primary advantage is that it is not based on the location of the critical point since we do not know the analogue of the critical point in two variables.  It has the disadvantage that many distinct points have the same coding. On the other hand  for points in the intersection $W^s(\alpha)\cap W^u(\alpha)$, the set of points with the same code is finite and this fact will be very useful for us. For example the two points in Figure~0.1 have the same coding sequence and, if we consider all the points corresponding to this coding sequence we obtain 4 points, which are marked ``$\circ$'' in Figure~0.2.  The degeneration of the horseshoe is shown schematically in Figure~0.3.  We note that the ``right-left'' coding with respect to the arc of $W^s(\beta)$ remains stable  throughout the bifurcation.

The maps we will work with are of the form $f_{c,\delta}(x,y) = (c +\delta y - x^2,-x)$.  
The values $(c,\delta)$ are taken from parameter region ${\cal W}_*$ for which we have a system of crossed mappings as described below.   Such parameter regions were constructed in [BSii] and especially in [BSh, \S2], where they are shown to contain a large part of the boundary of the real horseshoe locus (cf.\ [BSh, Figures 6 and 9]).  Since the parameter values themselves to not play any role in this paper, we will drop them from the notation and simply say that $f$ belongs to the region ${\cal W}_*$.

Our main theorem shows how the horseshoe and maximal entropy properties can be located in parameter space by using this coding.

\proclaim Theorem 1.  Let $f$ be an orientation-preserving real H\'enon map in the region $\cal W_*$.
Consider the collection of (real)  points with coding sequence $\overline 0 101\overline 0$. There are at most 4 such points, and:
\item{(i)}  If there are exactly 4 such points, then $f$ is a hyperbolic horseshoe, by which we mean that it is hyperbolic and conjugate to the full 2-shift.
\item{(ii)} If there are 3 such points, then $f$ has a quadratic tangency but entropy $\log 2$.
\item{(iii)} If there are less than 3 such points, then $f$ has entropy less than $\log 2$.

In the case of ``the first tangency,'' we are able to describe the dynamics more precisely:
\proclaim Theorem 2.  If $f$ is as in case (ii) above, then $f$ is topologically conjugate to the shift map on  $\Sigma_2/\sim$, which is the quotient  of the full shift on two symbols $\{a,b\}$, modulo the identification $\overline a bab\overline a\sim \overline abbb\overline a$.

We will make the situation of Figures 0.2 and 0.3 rigorous by passing to the complex domain.  
The idea of using complex techniques for the study of the real H\'enon map goes back to Hubbard. In particular, [HO] proves the existence of horseshoes using complex techniques.
In this paper we are interested in applying complex techniques to the analysis of
the real H\'enon family and particularly in locating the horseshoe region and its boundary and analyzing the topological dynamics of maps in the boundary. 

One of our long term goals is to understand how the beautiful theory of 1-dimensional complex dynamics can be applied in dimension 2.  Our method in this paper is related to the so-called ``puzzle'' corresponding to the period 2-wake in the parameter space of quadratic polynomials, which is constructed by using the external rays landing at the $\beta$ fixed point.  

The construction of wakes is based on the theory of external rays and the question of which collections of external rays land at given periodic points (see [M]).  The regions in parameter space for which a given landing pattern occurs are called wakes. For the parameters lying in a given wake the collection of rays landing at a periodic point can be used to partition dynamical space and produce a dynamical coding of points in this space.
The region in parameter space where our coding is defined is an analog of the period 2 wake in the parameter space of $z\mapsto z^2+c$. This is the region where the external rays of angles 1/3 and 2/3 land at the noncontracting fixed point. These rays divide complex dynamical space into two regions and give us a coding of points by strings of 0's and 1's.   When we work in the complex domain, this coding allows us to define specific local pieces of $W^s(\alpha)$ (or $W^u(\alpha)$) within each box. These local stable/unstable pieces are vertical (or horizontal) disks inside their respective boxes and they vary continuously with parameters. An interesting discussion of wakes and its relation to quadratic H\'enon maps is given in the thesis of Lipa [L].

\bigskip\noindent{\bf 1.  Definition and Properties of Crossed Maps. }
Let us describe the larger context of our work.  We consider an open set $U\subset{\bf C}^2$ and a holomorphic map $f:U\to {\bf C}^2$, which we consider as a partially defined dynamical system.  That is, we consider the dynamical system $f:K_\infty\to K_\infty$, where we define
$$K_\infty:=\{z\in U: f^nz\in U\ \forall n\in{\bf Z}\}. \eqno(1.1)$$
In the cases we consider in detail, $f$ will be a polynomial (and thus entire) mapping, and $K_\infty$ will coincide with the set of points with bounded orbits.  One useful tool for this study is the crossed mapping, which was introduced by Hubbard and Oberste-Vorth [HO].  If a crossed map from a domain to itself has degree 1, it exhibits hyperbolic behavior.  If it has higher degree, it is ``H\'enon-like''; such maps were studied by Dujardin [Du].  Our strategy in this paper is to replace $U$ by a covering $\{B_i\}$ and to replace  $f$ by a system of crossed mappings $f_{i,j}$ from $B_i$ to $B_j$.  Such an approach was adopted earlier in [BSii] and [DDS].  In the 3-box system, which we discuss below, one of the crossed maps will have degree 2, while all the others will have degree 1.

\medskip
\epsfysize=1.3in
\centerline{ \epsfbox{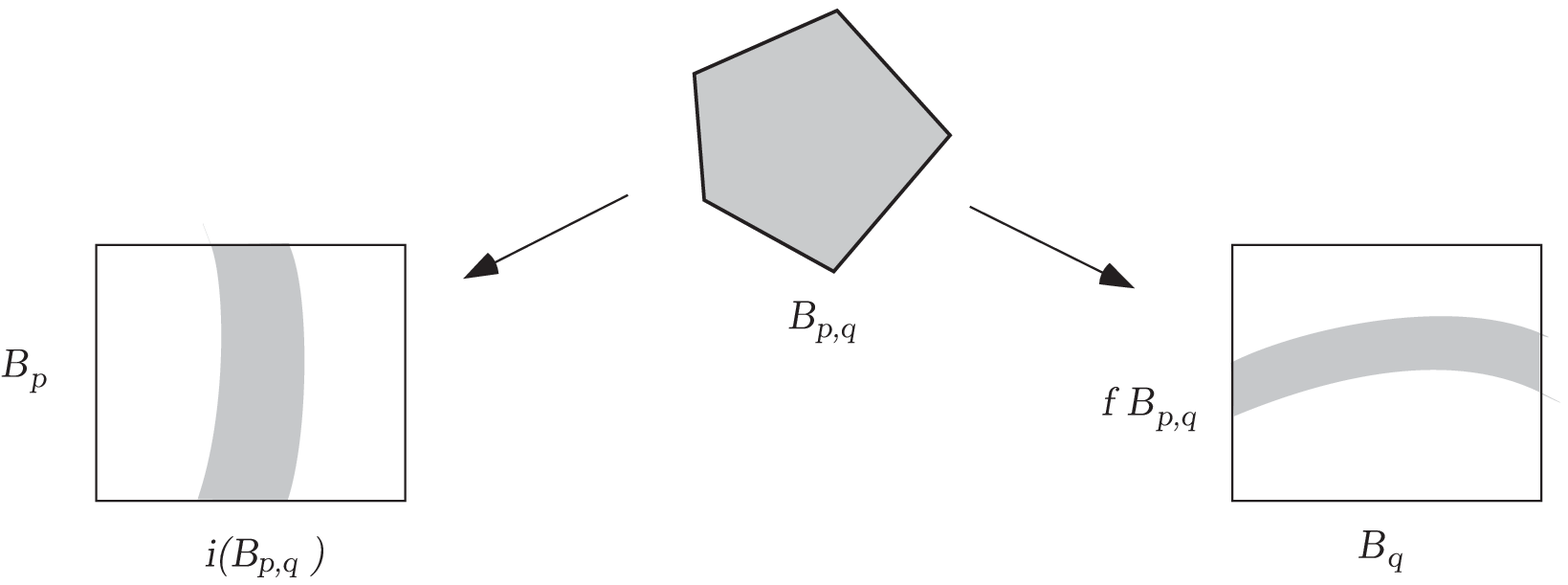}  }
\centerline{Figure 1.1  A Crossed Mapping}
\medskip

We observe that there will be some overlap between the definitions and notations presented in this Section and those in [IS].  Let $B_{p}$ and $B_{q}$ be bidisks.  For $j=1,2$, let $\pi_{j}:\Delta\times\Delta\to\Delta$ be the projection onto the $j$th coordinate.  We consider a complex manifold $B_{p,q}$ and a pair of holomorphic mappings
$$\iota=\iota_{p,q}:B_{p,q}\to B_{p}$$
$$f=f_{p,q}:B_{p,q}\to B_{q}.$$
We say that $(\iota, f,B_{p,q},B_{p},B_{q})$ is a {\it crossed mapping from $B_{p}$ to $B_{q}$} if
\item{(i)} $\iota B_{p,q}$ is holomorphically convex in $B_{p}$, and $f  B_{p,q}$ is holomorphically convex in $B_{q}$.

\item{(ii)} $((\pi_{2}\circ\iota),(\pi_{1}\circ f)):B_{p,q}\to\Delta\times\Delta$ is proper.

\item{(iii)}  There exists $\kappa<1$ such that $\pi_{1}(\iota B_{p,q})\subset\{|z|<\kappa\}$ and $\pi_{2}(f B_{p,q})\subset\{|z|<\kappa\}$.

\item{(iv)}  $\iota$ is an injection.

We will sometimes find it convenient to denote the crossed mapping simply by $f_{p,q}$.
We may define the composition of crossed mappings $f_{1,2}$ and $f_{2,3}$ as follows.  We set 
$$B_{1,3}:=\{(z,w)\in B_{1,2}\times B_{2,3}: f_{1,2}z=\iota_{2,3}w\},$$ 
and we define $\iota_{1,3}:=\iota_{1,2}|_{B_{1,3}}$, $f_{1,3}:=f_{2,3}|_{B_{1,3}}$.  The following is evident:

\proclaim Proposition 1.1.  The composition $f_{2,3}\circ f_{1,2}:=f_{1,3}$ is a crossed map from $B_1$ to $B_3$.

Let $\cB=\{B_{1},\dots,B_{N}\}$ be a family of bidisks, and let $\cG$ be a directed graph with vertices labeled $\{1,\dots,N\}$.  We say that a pair $(i,j)$ is {\it admissible} if there is an edge of $\cG$ going from $i$ to $j$.  We denote admissibility by $(i,j)\in\cG$.  We say that $(\cB,\cG)$ is a {\it system of crossed mappings} if there is a crossed mapping $f_{i,j}$ from $B_i$ to $B_j$ for each $(i,j)\in\cG$.   We can also deal with the situation where the graph $\cG$ has multiple edges.  In that case we use the edges of $\cG$ to label orbits.

We say that a sequence $s$ (finite or infinite) is {\it admissible} if it is compatible with $\cG$.  For instance, $s=(s_{n},s_{n+1},\dots,s_{m})$ is admissible if $(s_{j},s_{j+1})\in\cG$ for all $n\le j\le m-1$.  If $s$ is a finite, admissible sequence, then by Proposition 1.1 there is a unique crossed mapping $f_{s}:=f_{s_{m-1},s_{m}}\circ\cdots\circ f_{s_{n},s_{n+1}}$ from $B_{s_{n}}$ to $B_{s_{m}}$.  We will also refer to an admissible sequence as an {\it index}.  We let $\Sigma_{\cG}$ denote the set of bi-infinite admissible sequences.  We let $\Sigma_{\cG}[n,m]$ denote the admissible sequences for the interval $n\le j\le m$.

An {\it orbit} is a pair of sequences $(s,p_{s})$, where $s$ is an index, and $p_{s}=(p_{s_{n}},\dots,p_{s_{m}})$ is a sequence of points such that $p_{s_{j}}\in B_{s_{j}}$, and $f_{s_{j},s_{j+1}}p_{s_{j}} = p_{s_{j+1}}$ for all $n\le j\le m-1$.  We let $\Lambda$ denote the set of bi-infinite orbits, and we let $\Lambda(n,m)$ denote the finite orbits for which $j$ runs from $n$ to $m$.  There is a dynamical system on $\Lambda$ given by the shift map; and the projection $\pi:\Lambda\to\Sigma_{\cG}$ gives a semi-conjugacy to the shift on $\Sigma_{\cG}$.

Consider a holomorphic imbedding $\phi:\Delta\to B_{j}$.  We say that $\cD=\phi(\Delta)$ is a {\it horizontal disk} in $B_{j}$ if $\pi_{1}:\cD\to\Delta$ is proper.   In this case $\pi_{1}\circ\phi:\cD\to\Delta$ has a well-defined mapping degree $\mu$, and we say that $\mu$ is the {\it degree} of $\cD$.  We denote this as $\d(\cD)$.   Similarly, we define $\cD$ to be a vertical disk in $B_{j}$ if $\pi_{2}\circ\phi$ is proper, and we define the degree of a vertical disk to be the degree of this map.  A {\it horizontal multi-disk} is a union of horizontal disks.
Let $\cH(B_{j})$ denote the set of horizontal disks inside $B_{j}$.  We let $\cH_{\mu}(B_{j})$ denote the set of horizontal multi-disks in $B_{j}$ whose total degrees add to $\mu$.
\bigskip
\epsfysize=.9in
\centerline{ \epsfbox{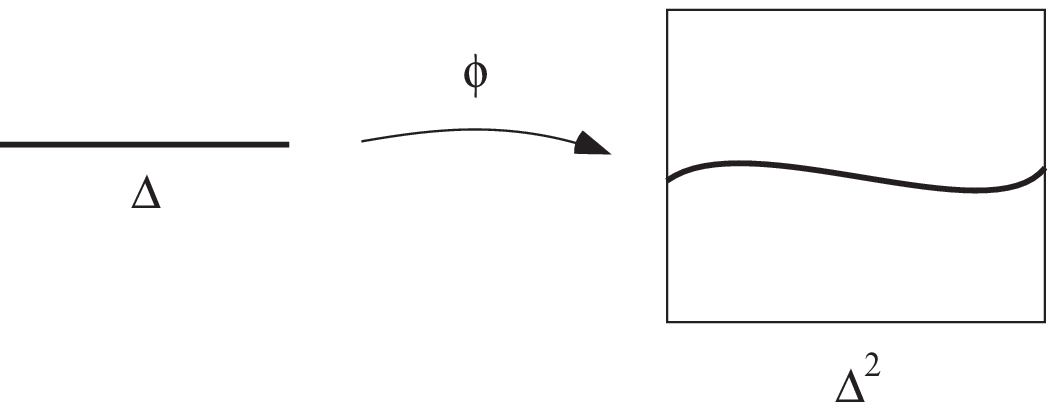}  }
\centerline{Figure 1.2.  A Horizontal Disk}

\proclaim Proposition 1.2.  Let $\cD_{h}$ and $\cD_{v}$ be horizontal and vertical multi-disks in $B_{j}$.  Then the number of intersection points (counted with multiplicity) is
$$\#(\cD_{h}\cap\cD_{v})=\d(\cD_{h})\d(\cD_{v}).$$

We note that Proposition 1.2 and the other results in the rest of this section are parallel to results in [HO], and we omit the proofs to avoid unnecessary duplication.

\proclaim Proposition 1.3.  Let $(f_{j,k},B_{j},B_{k})$ be a crossed map.  Then there is a number $\d(f_{j,k})$ such that if $\mu'=\d(f_{j,k})\mu$, then
$$f_{j,k}:\cH_{\mu}(B_{j})\to\cH_{\mu'}(B_{k})$$
and 
$$f^{-1}_{j,k}:\cV_{\mu}(B_{k})\to\cV_{\mu'}(B_{j}).$$

We define the {\it degree of a transition} $(i,j)\in\cG$ to be $\delta(i,j):=\d(f_{i,j})$.  If $s=(s_{n},\dots,s_{m})$ is a finite index, we define the degree 
$$\delta(s)=\delta(s_{n},s_{n+1})\delta(s_{n+1},s_{n+2})\dots\delta(s_{m-1},s_{m}).$$
For infinite $s$, we define $\delta(s)$ similarly if all but finitely many transitions have degree 1.
Thus {\it degree} is defined for disks, maps, and indices.

\proclaim Proposition 1.4.  If $s$ is a finite index, then $\d(f_{s})=\delta(s)$.

If $s=(s_{n},\dots,s_{m})$, $n\le0\le m$, is a finite index, we define
$$B(s):=\{z\in B_{s_{0}}:\exists{\rm\ orbit\ }(s,(p_{s_{n}},\dots,p_{s_{m}})){\rm\ with\ }p_{s_{0}}=z\}.$$
If $s$ starts at $j=0$, i.e. $s=(s_{0},\dots,s_{n})$, then $B(s)$ is the domain of definition of $f_{s}$.  If $s$ is infinite, we set $B(s)=\bigcap_{s'}B(s')$, where the intersection is taken over all finite sub-indices $s'$ of $s$.

\proclaim Theorem 1.5.  Let $s=(s_{0},s_{1},s_{2},\dots)$ be a semi-infinite index of degree 1.  Then $B(s)$ is a horizontal disk inside $B_{s_{0}}$.  Further, $B(s)$ is an unstable manifold: any two points of $B(s)$ approach each other in backward time.

\proclaim Theorem 1.6.  Let $s=(s_{0},s_{1},s_{2},\dots)$ be a semi-infinite index of finite degree.  Then $B(s)$ is a horizontal multi-disk inside $B_{s_{0}}$.  Further, $B(s)$ is a piece of an unstable manifold.

We say that a finite index $s=(s_{0},\dots,s_{m})$ is {\it periodic} if $s_{0}=s_{m}$.   A periodic index represents a closed path in $\cG$.  For each periodic index $s$ there is a crossed map $(f_{s},B(s),B(s))$.  

\proclaim Theorem 1.7.  If $s$ is a periodic index with degree 1, then the induced crossed mapping $f_{s}$ of $B(s)$ has a unique fixed point.  This fixed point is a saddle point.

\bigskip\noindent{\bf 2.  Codings of Orbits. }  
Let an open set  $U\subset\C^{2}$ and a holomorphic map $f:U\to\C^{2}$ be given.    We say that a system of crossed mappings $(\cB,\cG)$ is a {\it realization } of the partially defined system  $(f,U)$ if there is a family of holomorphic imbeddings $\phi_{j}:B_{j}\to U$, $1\le j\le N$  with the property that $U=\bigcup_j \phi_jB_j$, and for each $(j,k)\in\cG$,  $\phi_{k}\circ f_{j,k}=f\circ\phi_{j}$ holds on the common domain of the two maps.  
If $\phi$ is a realization as above, then there is a semiconjugacy $\phi:\Lambda\to K_{\infty}$ defined by $\phi(s,p_{s})=\phi_{s_{0}}(p_{s_{0}})$, where $K_\infty$ is defined in (1.1). 

When there is no danger of confusion, we will drop the $\phi_j$, writing simply $B_j$ for $\phi_jB_j$.   We use the notation
$$K_\infty=\{z\in U: f^nz\in U\ \forall\ n\in{\bf Z}\}.$$

We use the notation $\dot B_j$ to denote the smaller bidisk $\{|z|<\kappa\}\times\Delta$.  As before, we also write $\dot B_j$ to denote $\phi_j\dot B_j$ and $\dot U:=\bigcup \dot B_j$.  With obvious notation, we have crossed mappings $\dot f_{i,j}$ from $\dot B_i$ to $\dot B_j$.    We may choose $\eta>0$ such that
$$\dot U{\rm\ contains\ an\ }\eta{\rm\, neighborhood\ of\ }K_\infty, {\rm \ and\ }U\supset \{z\in fU:{\rm\ dist}(z,\dot U)<\eta\}.\eqno(2.1)$$
We will also want the system $(\dot\cB,\cG)$ to satisfy the following two properties:
$$\eqalign{ &\forall \ i {\rm\ and\ }z\in \dot B_i\cap f\dot U, \ {\rm\exists\ }\alpha{\rm\ such\ that\ }    \cr
&z\in \dot f_{\alpha,i}\dot B_\alpha{\rm\ and\  dist}(z,\partial(\dot  f_{\alpha,i}\dot B_\alpha)\cap\partial  \dot B_i)>\eta.\cr} \eqno(2.2)$$ 
$$\forall\ \alpha',\alpha'' ,i \ {\rm with}\ ( \alpha',i),(\alpha'',i)\in\cG, \ {\rm and }\ \alpha'\ne\alpha'',\ 
{\rm dist}(\dot f_{\alpha',i}\dot B_{\alpha'},  \dot f_{\alpha'',i}\dot B_{\alpha''})>\eta.  \eqno(2.3)$$ 
We note that condition (2.2) says, in some sense, that the system $(\dot \cB,\cG)$ has ``enough'' maps to mirror the dynamics of $f$, and condition (2.3) says that the forward images of the $B_j$'s are vertically separated.  

We say that a (finite or infinite) sequence  $(z_n)$  in $U$ is an $\epsilon$ orbit  (or $\epsilon$ pseudo-orbit) if ${\rm dist}(z_{j+1},fz_j)<\epsilon$ for all $j$.  Since the sets $B_j$ which make up $U=\bigcup_j B_j$ can have nonempty pairwise intersection, the condition that $z\in B_j$ does not uniquely determine $j$.  However, if (2.1--3) are satisfied, then we can define a coding which is consistent with the itinerary of an $\epsilon$ orbit.
\proclaim Theorem 2.1.  For a realization satisfying (2.1--3), there is an $\epsilon>0$ with the following property.  If $z_0,\dots,z_n\in U$ is an $\epsilon$ pseudo-orbit for $f$, and if $fz_n\in U$, then there is an admissible sequence $j_0,\dots,j_n$ such that $z_i\in \phi_{j_i}B_{j_i}$ for each $0\le i\le n$.  That is, the sequence of pairs $(z_i,j_i)$ is an $\epsilon$ pseudo-orbit for the system of crossed mappings.  Further, if $z_n\in f_{j,\alpha}B_j$ for some $\alpha$, then there is a unique sequence $j_0,\dots,j_n$ which ends in $\alpha$.   Similarly, if $(z_n)$ is an infinite $\epsilon$ pseudo-orbit in $U$, then it is coded by an infinite sequence $(j_n)$.  In either case, there are at most $N$ such sequences of symbols.

\noindent{\it Proof. }  We take $\epsilon<\eta/2$.  First, since $z_n\in U=\bigcup B_i$, we may choose an index $i_n$ such that $z_n\in B_{i_n}$.  Now suppose that $n>1$.  Since $z_{n-1},z_n\in U$, it follows from (2.1) that $fz_{n-1}\in\tilde B_{i_n}\cap f\tilde U$.  Now by (2.2) there is an index $i_{n-1}$ such that $(i_{n-1},i_n)$ is admissible, and $z_{n-1}\in\tilde B_{i_{n-1}}$.  Further, by (2.3), the index $i_{n-1}$ is uniquely determined by the choice of $i_n$.  For $\alpha\ne i_{n-1}$, and if $(\alpha,i_n)$ is admissible, then the separation between $f_{\alpha,i_n}\tilde B_\alpha$ and $f_{i_{n-1},i_n}\tilde B_{i_{n-1}}$ is greater than $\eta$.  Continuing backwards in this fashion, we obtain the admissible sequence $i_0,\dots,i_n$.  Thus the number of possible sequences $(i_j)$ is no greater than the number of possible choices  of $i_n$, which in turn is no greater than $N$.

Now suppose the $\epsilon$ orbit $(z_n)$ is infinite.  For each number $k$, we may choose an index $i_k^{(k)}$ such that $z_k\in B_{i_k^{(k)}}$.  As above, we can work our way to the left, and thus construct an admissible sequence $s^{(k)}=\dots\, i_j^{(k)}\dots  \,i_k^{(k)}$.  Since there are at most $N$ possible sequences at each stage $k$, we see that there will be an infinite sequence of values $k_m\to\infty$ such that the sequences $s^{(k_m)}$ extend each other.  This gives the desired infinite sequence.  \qed

We we may restate part of the Theorem as follows: 
\proclaim Corollary 2.2.  The map $\phi:\Lambda\to K_\infty$ is a surjection, and for $p\in K_\infty$, $\#\phi^{-1}p\le N$.

The proof of Theorem~2.1 only used properties of the $\tilde B_i$ which remain valid under small perturbations of both the pseudo-orbit and the mapping.  Thus we have continuity for the coding process:
\proclaim Theorem 2.3.  Let $\epsilon>0$, and let $z_0,\dots,z_n$, and $z_0',\dots,z_n'$ be two $\epsilon$ pseudo-orbits 
as in the previous Theorem.  There is a $\delta>0$ such that if any two such pseudo-orbits are $\delta$ close, then they have the same codings.

If a system of boxes works for a map $f$, then it will work for a small perturbation of $f$.  Thus we have:
\proclaim Theorem 2.4.  If (2.1--3) holds for $f$, and if $f'$ is sufficiently close to $f$ on $U'$, then the coding of a the orbit of a point with respect to $f$ is the same as its coding with respect to~$f'$.

\bigskip\noindent{\bf 3.  The 3-Box System. }
Let us consider a system of crossed mappings $(\cB,\cG)$, which consists of 3 bidisks  $\cB=\{B_{0},B_{1},B_{2}\}$ and the graph  $\cG$ pictured in Figure 3.1.  The degrees of the crossed mappings are required to satisfy $\delta(1,2)=2$, and $\delta(i,j)=1$ for all other $(i,j)\in\cG$.   Let $f:U\to {\bf C}^2$ be a holomorphic mapping associated to $\cB$ as in \S2.  Finally, we require that $B_{(12)}:=\phi_1B_1\cap \phi_2 B_2$ is a bidisk with respect to the product structures of both $B_1$ and $B_2$ and has the form $\Delta_{(12)}\times \Delta$, i.e., it extends to full height in $B_1$ and $B_2$.  The map $f$ then defines a degree 1 crossed mapping of $B_{(12)}$ to itself.  If all of this holds, and if conditions (2.1--3) are satisfied, we call $(\cB,\cG)$ a 3-box system for the map $f:U\to{\bf C}^2$.
\bigskip
\epsfysize=.6in
\centerline{ \epsfbox{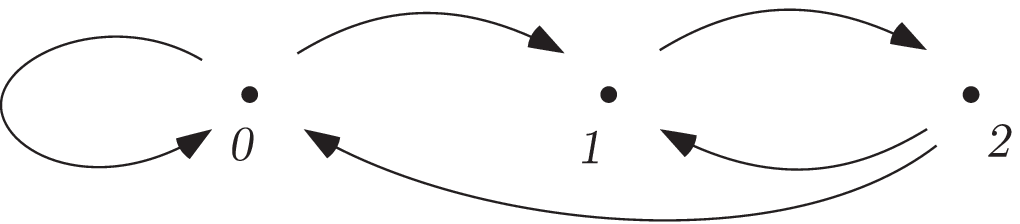}  }
\centerline{Figure 3.1  Graph $\cG$ of the 3-box system.}
\bigskip
\epsfysize=.7in
\centerline{ \epsfbox{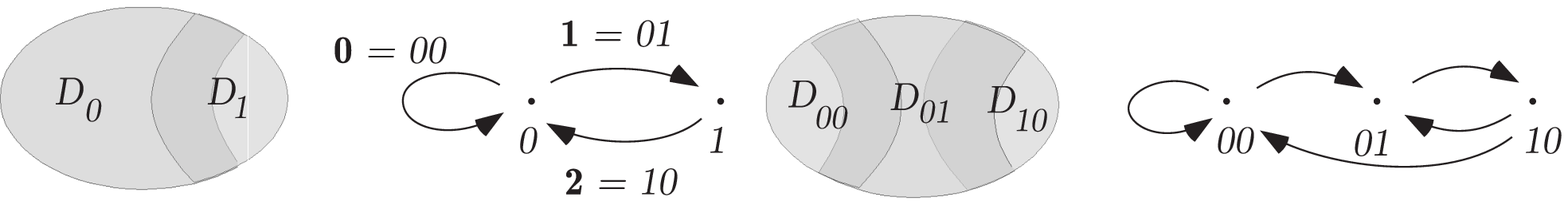}  }
\centerline{Figure 3.2  Box systems: Graph $\cH$ on the left, and 2-orbits on right.}
\medskip
In [BSii] and [BSh] 3-box systems were constructed by first making a 1-D system and then taking a product with $\Delta$.  This construction proceeds by starting with a map $p:{\bf C}\to{\bf C}$ and finding a 2-box system $(\cD,\cH)$, $D_0,D_1\subset{\bf C}$, as on the left hand side of Figure~3.2.  Shading denotes overlap of boxes.  The orbits for the system $(\cD,\cH)$ can be coded in terms of a 01 coding corresponding to the graph $\cH$.  Then we  pass to the orbits of length 2 by setting $D_{j,k}:=\{z\in D_j:p(z)\in D_k\}$, and obtain a 3-box system as shown on the right hand side of Figure~3.2.  We see that the graph on the right hand side of Figure~3.2 is obtained as the edge coding of the graph $\cH$.  The codings in terms of the two graphs in Figure~3.2 are essentially the same, and we may pass between the $\cG$- and the $\cH$-codings by the correspondence ${\bf 0}\leftrightarrow{00}$, ${\bf 1}\leftrightarrow{01}$, ${\bf 2}\leftrightarrow{10}$.  (While we are discussing both the $\cG$- and the $\cH$-codings, we use bold face to identify the $\cG$-coding.)  While the $\cH$-coding is more natural from certain points of view, we find it convenient to have 3 symbols, 012, to identify the three boxes.  The index  $\overline 0101\overline 0$ in the statement of the Main Theorem is given in the $\cH$-coding; the corresponding $\cG$-coding, which we will use extensively in the sequel, represents this index as $\overline{\bf 0}{\bf 1212}\overline{\bf 0}$.

\bigskip

\epsfysize=.6in
\centerline{ \epsfbox{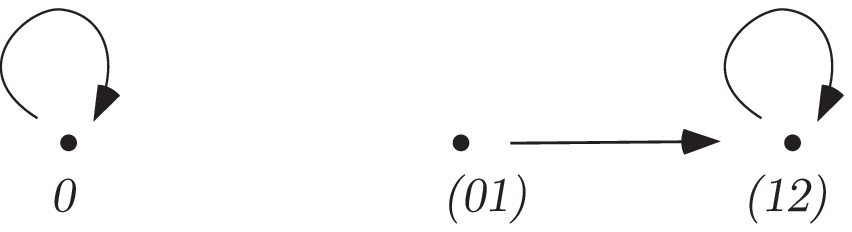}  }

\centerline{Figure 3.3.  Degree One Subsystems.} 

\medskip
We observe that a 3-box system $(\cB,\cG)$ contains the degree one subsystem given by the crossed map from $B_0$ to itself.  In addition, we set $B_{(01)}:=\phi_0B_0\cap \phi_1B_1$, and we have the degree 1 system from $B_{(01)}$ to $B_{(12)}$, and the crossed map of $B_{(12)}$ to itself, as pictured in Figure 3.3.   It follows from \S1 that there is a unique saddle point $p_0\in B_0$, and the local stable/unstable manifolds, $W^s_{loc}(p_0)=B({\bf 0}\overline 0)$ and $W^u_{loc}(p_0) =B(\overline 0{\bf 0})$ are vertical/horizontal disks of degree 1 in $B_0$.  Here and in the sequel we will change our notational convention slightly and use bold face type to indicate which index is in the zero-th position. 

There is a unique saddle $p_{(12)}\in B_{(12)}$, and the (unique) stable manifold of the system $(f_{(12)},B_{(12)})$ is a degree one vertical disk, shown as the arc $\sigma$ in Figure 5.3.  It is evident that $\sigma$ is contained in both $B^r({\bf 1}\overline{21})$ and $B^r({\bf 2}\overline{12})$.  The curve $\sigma$  is vertical and degree 1 in both $B_1$ and $B_2$.  Thus $f_{0,1}^{-1}\sigma\subset B_0$ is vertical of degree 1, and $f^{-1}_{1,2}\sigma\subset B_1$ is vertical of degree 2.  Since $\sigma$ is a stable manifold,  there is a vertical curve $\sigma'\subset B_0\cap B_1$ such that $f_{1,2}^{-1}\sigma=\sigma\cup \sigma'$, and $f^{-1}_{0,1}\sigma=\sigma'$.  The corresponding boxes are $B^r({\bf 0}\overline{12})=\sigma'$, $B^r({\bf 1}\overline{21})=\sigma\cup\sigma'$, and $B^r({\bf 2}\overline{12})=\sigma$.

Now we ask which points  $p\in K_\infty$ can have non-unique codes for the 3-box system.  If $f^np\notin B_{(01)}\cup B_{(12)}$, then by Theorem 2.1, the symbols $s_j$ which code the orbit are uniquely determined for $j\le n$.  Thus there must be an $n_0$ such that $f^np\in B_{(01)}\cup B_{(12)}$ for  $n\ge n_0$.   Thus for $p$ to have a non-unique coding, its orbit must enter the degree 1 subsystem and remain inside it for all forward time.  We summarize these remarks on non-uniqueness of coding as follows:
\proclaim Proposition 3.1.  If $p\in K_\infty- W^s(p_{(12)})$ then $\#\phi^{-1}p=1$; and if $p\in K_\infty\cap W^s(p_{(12)})$, then $\#\phi^{-1}p=2$.

Now let us give two ways of shrinking the boxes $B_j$.  
We may shrink them in forward time by setting
$$B_0^+:= \bigcup_{j\ge0} B({\bf 0}0^j1), \ \ B_1^+:=\bigcup_{j\ge0} B({\bf 1}(21)^j20) \ \ {\rm and}\ \  B_2^+:=\bigcup_{j\ge0}B({\bf 2}(12)^j0),$$
where we use the notation $0^j=0\cdots 0$ and $(12)^j=12\cdots 12$.  We observe that for any $n$,
$$K_\infty\cap B_0\subset B_0^+\cup B({\bf 0}0^n),\ \  K_\infty\cap B_1\subset B_1^+\cup B{\bf 1}(21)^n), \ \ {\rm and}\ \  K_\infty\cap B_2\subset B_2^+\cup B({\bf 2}(12)^n).$$
 $B({\bf 0}0^n)$ is an open set which shrinks down to the local stable manifold $B({\bf 0}\overline 0)$ as $n\to\infty$.  Similarly, $ B({\bf 1}(21)^n)$ and $B({\bf 2}(12)^n)$ are neighborhoods which shrink down to the local stable manifold $\sigma$ as $n\to\infty$.  $B({\bf 0}1)$ is equivalent to  a product of two disks with the base disk lying inside the base of $B_0$ as shown on the left hand side of Figure 3.4.  Further, $B({\bf 0}0^n1)=f_{0,0}^nB({\bf 0}1)$.  Since $p_0$ is a saddle point, we may apply the Lambda Lemma to see that $f_{0,0}^{-n}B({\bf 0}1)$ is essentially the product of two disks, and one base disk is mapped to the next by a map which is approximately conformal. 
 In one variable complex dynamics John conditions (also called carrot conditions) on Julia sets are related to expanding properties of the map. Though we make no further  use of this fact in this paper we observe that we get some John type conditions.  In particular the complement $B_0-B^+_0$, and thus $B_0-K_\infty$,  satisfies a John condition at every point of the local stable manifold $B({\bf 0}\overline0)$. 
 In fact, we have a very uniform family of carrots since  $B_0-K_\infty$ contains a carrot-wedge that is essentially the product of a carrot (in the base disk and landing at $p_0$) and a (vertical) disk.  The situation in the base disk is shown on the right hand side of Figure 3.4: clearly, there is a carrot landing at $p_0$, which is disjoint from the union of the disks representing the bases of $B({\bf0}0^j1)$.
\medskip

\epsfysize=1.0in
\centerline{ \epsfbox{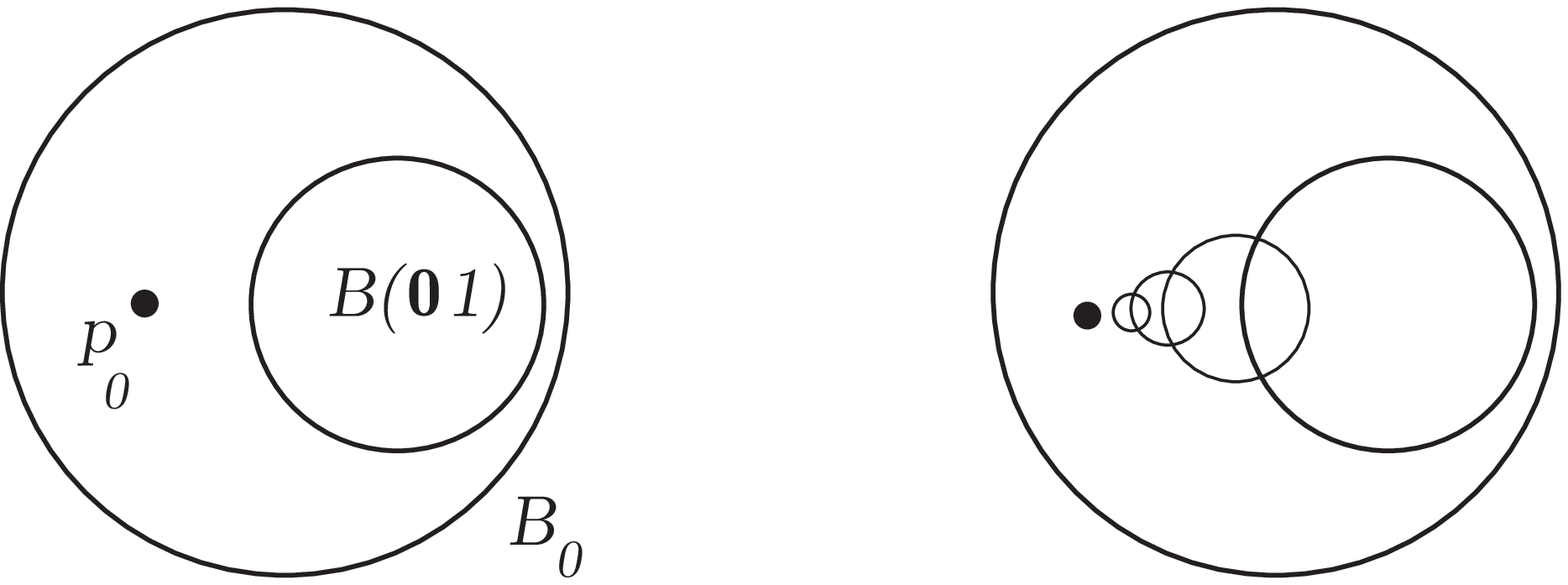}  }

\centerline{Figure 3.4.  John condition in the base disk.} 
\medskip
Similarly, we define $B^-_0:=\bigcup_{j\ge0} B(20^j{\bf 0})$, and we see that $B_0-K_\infty$ satisfies a carrot condition at each point of the local unstable manifold $B(\overline0{\bf 0})$ of $p_0$.  Thus we have:
\proclaim Proposition 3.2.  Suppose that $f:U\to {\bf C}^2$ is realized by a 3-box system.  Then $U-K_\infty$ satisfies a John (carrot) condition at each point of $W^u(p_0)$, $W^s(p_0)$, and $W^s(p_{(12)})$.

\noindent{\bf 4.  Real, Crossed Mappings; the Disk Property.}
Let us define real, crossed mappings.  Let $\tau$ be an anti-holomorphic involution which respects the product structure of $B_{j}=\Delta\times\Delta$ for all $1\le j\le N$.    We say that a crossed mapping $(f_{j,k},B_{j},B_{k})$ is {\it real} if $f_{j,k}$ commutes with $\tau$.  Let $B_{j}^{r}=\{z\in B_{j}:\tau z=z\}$ be the real points of $B_{j}$.  Thus we may identify $B^{r}_{j}$ with the product $I\times I$ of real intervals.  We say that a complex disk $\cD$ in $B_{j}$ is {\it real} if $\tau\cD=\cD$.

\proclaim Theorem 4.1.  If $\cD$ is a real disk in $B_{j}$, then $\cD\cap B_{j}^{r}$ is a nonempty, connected, smooth arc.

The theorem is a consequence of the following.

\proclaim Proposition 4.2. The fixed point set of an orientation reversing diffeomorphism $\tau$ of a disk consists of a single arc.

\noindent{\it Proof of Theorem 4.1. } 
The restriction of $\tau$ to $\cD$ is an anti-conformal involution of the disk. The result now follows because the fixed points of $\tau$ are precisely the real points of $\cD$. \qed

\noindent{\it Proof of Proposition 4.2. }    By averaging we may assume that $\tau$ preserves some metric on the disk. By using the exponential map for that metric it follows that at a fixed point the involution $\tau$ is locally conjugate to the linear involution $D\tau$. A linear involution of $\R^2$ which reverses orientation has a one dimensional $+1$ eigenspace and a one dimensional $-1$ eigenspace. The $+1$ eigenspace of $D\tau$ is tangent to the fixed point set so the fixed point set is a one-manifold (perhaps with boundary). By choosing a triangulation compatible with the involution it is easy to see that the Euler characteristic of the fixed point set is congruent to the Euler characteristic of the disk mod 2 hence it is odd. Since a one-manifold is a union of circles, which have Euler characteristic 0, and arcs, which have Euler characteristic 1, we see that there must be at least one arc in the fixed point set. An arc in the disk divides the disk into two components. By looking at the action of $\tau$ near a fixed point we see that these components must be interchanged by $\tau$. In particular there can be no fixed points in the complement of this arc.   \qed

Note that if $\cD$ is a horizontal disk, then $\cD\cap B_{j}$ is a horizontal arc.  We define the intersection number $\#(\cD_{1}\cap\cD_{2}\cap B_{j}^{r})$ to be the number of intersection points of $\cD_{1}\cap\cD_{2}$ which are contained in $B_{j}^{r}$, and we count them with the multiplicity of the (complex) intersection between $\cD_{1}$ and $\cD_{2}$.  It is evident that
$$\#(\cD_{1}\cap\cD_{2})\ge\#(\cD_{1}\cap\cD_{2}\cap B_{j}^{r}).$$  
Since the complex intersection multiplicity of a point is always greater than or equal to one, equality holds if and only if $\cD_{1}\cap\cD_{2}\subset B_{j}^{r}$.

For comparison, we recall  the (oriented) real intersection number between the arcs $\cD_{1}\cap B_{j}^{r}$ and $\cD_{2}\cap B_{j}^{r}$  depends on the choice of orientation of each arc.  While the intersection number $\#(\cD_{1}\cap\cD_{2}\cap B_{j}^{r})$  may differ from the real intersection number, the two numbers are congruent modulo 2.

\bigskip\noindent{\bf 5.  The  3-Box System: The Real Case. }  Now let us describe a {\it real, 3-box system}.  This is a 3-box system $(\cB,\cG)$ with a real structure $\tau$, which is to say that $f$ preserves ${\bf R}^2$.  For convenience we will suppose,  in addition, that $f|{\bf R}^2$ preserves orientation.  While a parallel treatment for the non-orientable case is possible (see [BSii]),  we do not discuss it in this paper. 
 Let  $(i,j)\in\cG$ be an index with $\delta(i,j)=1$.  We know that $f_{i,j}$ maps (degree one) real, horizontal arcs in $B_{i}$ to (degree one) real, horizontal arcs in $B_{j}$.  We say that $f_{i,j}$ {\it preserves horizontal orientation} if it preserves the orientation of real, horizontal arcs.   For a real, 3-box system, we require that the horizontal orientation is preserved/reversed according to the $+/-$ signs in Figure 5.1.  Thus $f_{0,0}$ preserves horizontal orientation, and $f_{1,2}$ reverses it.  We require that the map
$f_{1,2}$ maps a degree one horizontal arc $\gamma$ in $B_{1}$ to $B_{2}$ according to one of the diagrams on the right hand side of Figure 5.2.  That is, $f_{1,2}\gamma$ is an arc of degree 2 in $B_{2}$, and we require that it intersect the left-hand boundary of $B_{2}$ in two points; and the second one is  below the first.

\smallskip
\epsfysize=.9in
\centerline{ \epsfbox{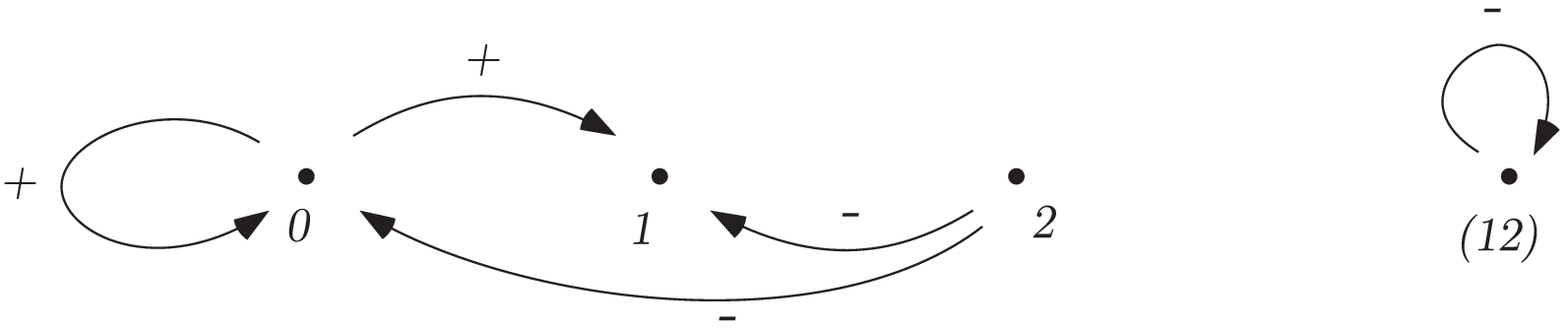}  }
\centerline{Figure 5.1: Preservation/reversal of horizontal orientation in the real, 3-box system.}
\medskip

Without loss of generality, we may assume that $B^r_0$ is arranged according to  Figure~5.3, i.e., $B^r({\bf 0}0)$ and $B^r({\bf 0}1)$ are (topologically) vertical strips
with $B^r({\bf 0}1)-B^r({\bf 0}0)$ to the right of $B^r({\bf 0}0)-B^r({\bf 0}1)$; and $B^r(0{\bf 0})$ and $B^r(2{\bf 0})$ are horizontal strips with $B^r(0{\bf 0})-B^r(2{\bf 0})$ above  $B^r(2{\bf 0})-B^r(0{\bf 0})$.  Thus we have the following:

\epsfysize=1.1in
\centerline{ \epsfbox{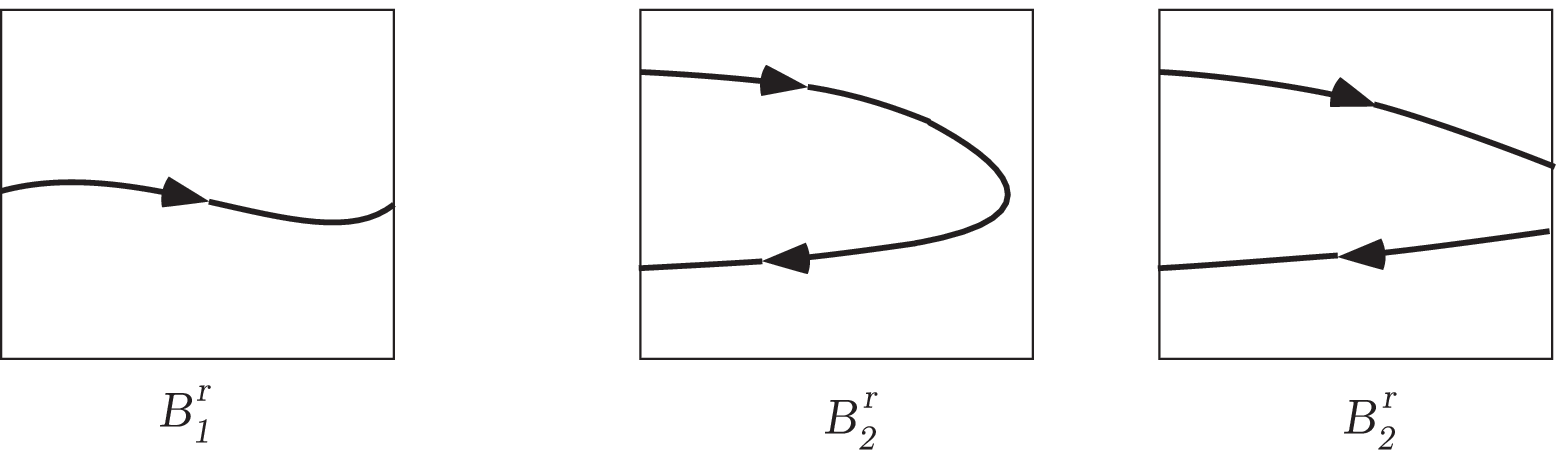}  }
\centerline{Figure 5.2: Geometry of the mapping $f_{1,2}$.}
\medskip
\epsfysize=1.3in
\centerline{ \epsfbox{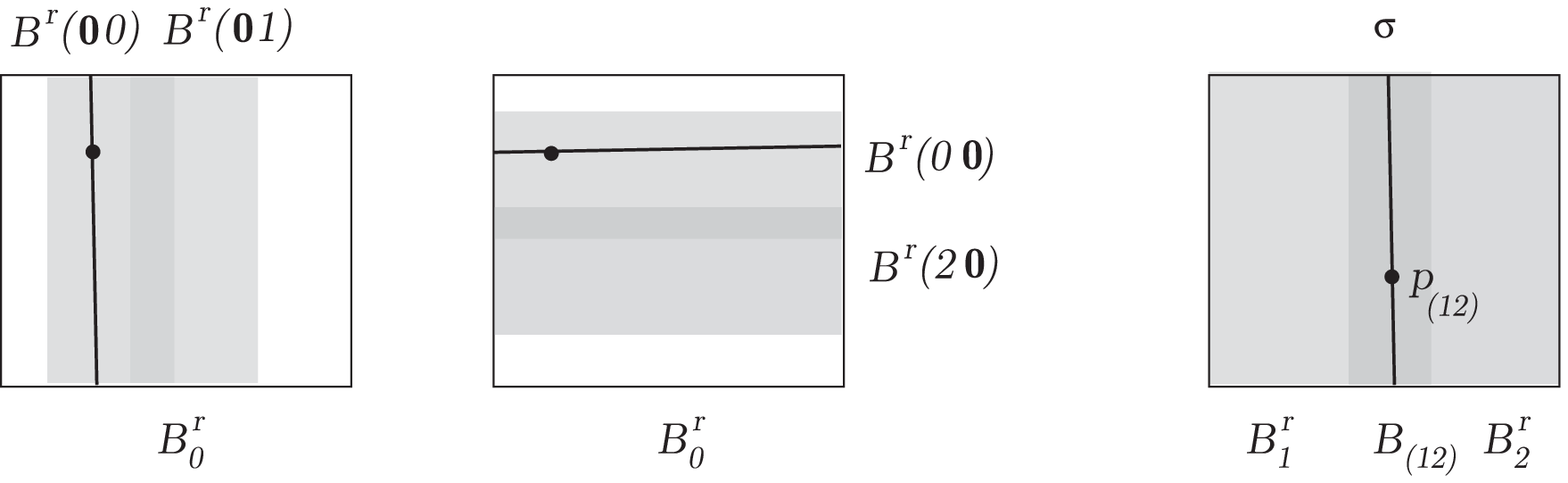}  }

\centerline{Figure 5.3.  Arrangement of $B_j^r$, $0\le j\le 2$.  }

\proclaim Lemma 5.1.  Every point of $B_0^r\cap K_\infty$ lies on or below the local unstable manifold $B^r(\overline 0{\bf 0})$ and to the right of $B^r({\bf 0}\overline 0)$.  No point of $B_0^r\cap K_\infty$ lies above $B^r(\overline 0{\bf 0})$, and  no point of $B_1^r\cap K_\infty$ lies above $B^r(\overline 0{\bf 1})$.  Every point of $B_2^r\cap K_\infty$ lies to the left of $B^r({\bf 2}\overline 0)$; and it  is enclosed by the loop of $B^r(\overline 0 1{\bf 2})$.

\noindent{\it Proof. }  Suppose that $p\in B_0^r\cap K_\infty$ is a point which is to the left of $B^r({\bf 0}\overline 0)$.  By the arrangement of the left hand side of Figure~5.3, we see that $p\notin B^r({\bf 0}1)$, which means that $fp\notin B_1^r$.  Thus we again have $fp\in B_0^r$.  The direction of horizontal orientation is preserved as in Figure~5.1.  Thus $fp$ continues to be to the left of $B^r({\bf 0}0)$.  Repeating this argument, we have $f^np\in B_0^r$ for all $n\ge0$, which means that $p$ lies in $B({\bf 0}\overline0)$, and not to the left of it.

Next we show that no point $p\in B_0^r\cap K_\infty$ lies above $B^r(\overline 0{\bf 0})$. By the central part of Figure~5.3, we see that $p\notin B^r(2{\bf 0})$, so $f^{-1}p\in B_0$.  By Figure~5.1, the map $f_{0,0}^{-1}$ must preserve horizontal orientation and, since $f^{-1}$ is orientation-preserving, $f_{0,0}^{-1}$ must also preserve vertical orientation.  Thus $f^{-1}p\in B_0^r$ and lies above $B(\overline 0{\bf 0})$.  Continuing this way we see that $f^{-n}p\in B_0^r$ for all $n\ge0$, which means that $p$ belongs to $B(\overline 0 {\bf 0})$ and is not above it. 

The assertion about $B^r_1\cap K_\infty$ now follows by mapping the result of the previous paragraph forward by $f_{0,1}$ and using the fact that $f$ must preserve vertical orientation.

The assertions about  $B^r_2\cap K_\infty$ being to the left of $B^r({\bf 2}\overline0)$ follow by mapping it back from $B_2$ to $B_0$.  And to see that $B_2^r\cap K_\infty$ lies inside   the loop of $B^r(\overline 0 1{\bf 2})$, we take the result of the previous paragraph and map it forward, keeping the geometry of Figure~5.2 in mind.  \qed

Now  let us consider the condition 
$$B^r(\overline 01{\bf 2}\overline0)\ne\emptyset. \eqno(*)$$
If  $(*)$ holds, the degree 2 horizontal curve $B^r(\overline 01{\bf 2})$ comes from the left and crosses the vertical line $B^r({\bf 2}\overline0)$.  In this case, $B(\overline01{\bf 2})$ contains sub-arcs $\tau'$ and $\tau''$ with the property that  $\tau'$ and $\tau''$ intersect $B^r({\bf 2}\overline0)$ exactly once and connect it to the left vertical side of $\partial B_2^r$.  If we map forward by $f$ we see that $f\tau'$ and $f\tau''$ will cross $B_1^r$ and connect $\sigma$ to the stable arc $B({\bf 0}\overline0)$.  Thus for $0\le j\le 2$ we may cut down the boxes $B^r_j$ to $B_j'$ as pictured in Figure~5.4.

\epsfysize=1.3in
\centerline{ \epsfbox{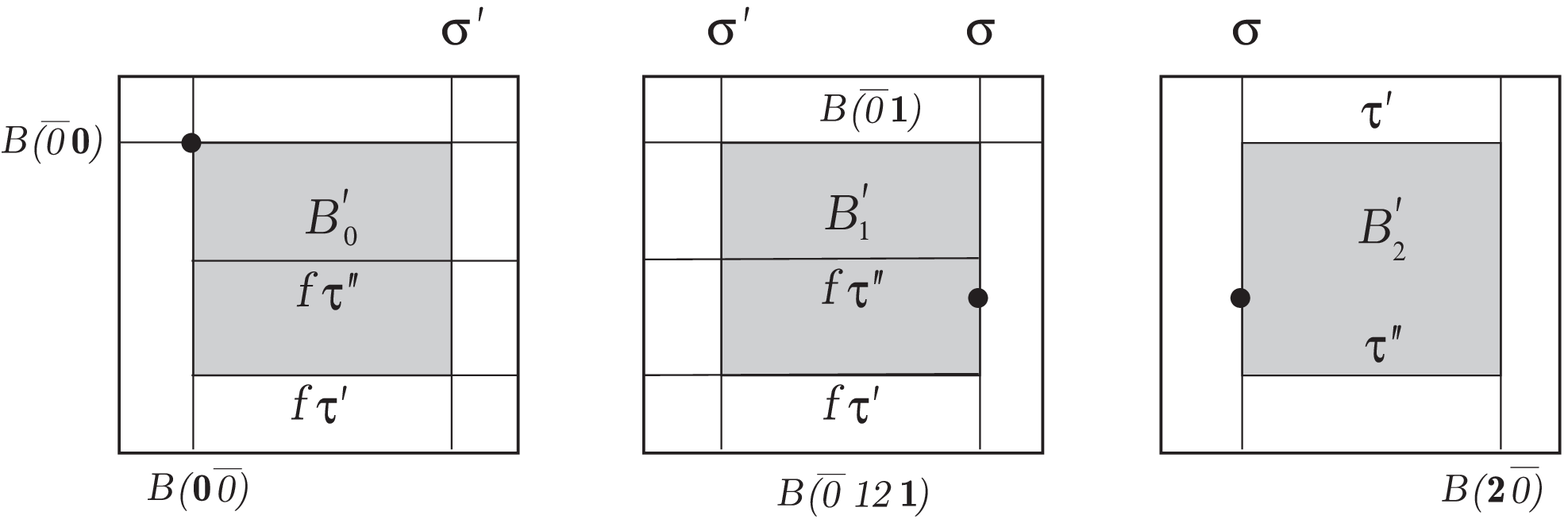}  }

\centerline{Figure 5.4.  Boxes $B_j'$, $0\le j\le 2$, are shaded.  }
\medskip

\proclaim Lemma 5.2.  If $(*)$ holds, and if $j=0$ or $j=1$, then no point of $B_j^r\cap K_\infty$ lies below $B(\overline0 12{\bf j})$.

\noindent{\it Proof. }  Since $(*)$ holds, we may construct the curves $\tau'$ and $\tau''$ as above.  By Lemma~5.1, $B_2^r\cap K_\infty$ lies below $\tau'$.  Since $f|B_2$ preserves total orientation and reverses horizontal orientation, we find, then, that $B_j^r\cap K_\infty$ must lie on or above $f\tau'$, which is the bottom portion of  $B(\overline0 12{\bf j})$. \qed

An immediate consequence of Lemmas 5.1 and 5.2 is the following:
\proclaim Lemma 5.3.  If $(*)$ holds, then $K_\infty\subset B_0'\cup B_1'\cup B_2'$.

Let us define $S:=B^r(\overline 0 121{\bf 2}\overline 0)$.  The degree of the index is $\delta(\overline 0 1212\overline 0)$ which is $4$, so the multiplicity of the complex box $B(\overline 0 1212\overline 0)$ is $4$.  Thus $\#S\le 4$.  

\proclaim Lemma 5.4.  If $\#S\ge 1$, then $(*)$ holds.

\noindent{\it Proof. }  By definition, $S\subset K_\infty\cap B_2^r$.  By Lemma 5.1, then, $S$ lies to the left of $B^r_2(\overline 01{\bf 2})$.  Further, $S\subset B^r({\bf 2}\overline0)$, so  $B^r(\overline 01{\bf 2})$ must intersect $B^r({\bf 2}\overline0)$, which means that $(*)$ holds.  \qed
\proclaim Lemma 5.5.  If $\#S\ge3$,  then  there are boxes $B_j^a$, $B_j^b$ as in Figure 5.5, with $B_j'\supset B_j^a\cup B_j^b$,   $B_j^a\cap B_j^b=\emptyset$ for $j=0,1$, and $\bigcup_j(B_j^a\cup B_j^b)\supset K_\infty$.   Further, if $\#S=3$, then $B_2^a\cap B_2^b$ is a single point, and if $\#S= 4$, then  $B_2^a\cap B_2^b=\emptyset$. 

\noindent{\it Proof. }  Let us start by showing that the construction of $B^{a/b}_j$, $0\le j\le 2$, indicated in Figure~5.5, is consistent.  We have seen already that condition $(*)$ implies that $B^r(\overline01{\bf 2})$, which has degree 2, contains arcs $\tau'$ and $\tau''$.  The arcs $f\tau'$ and $f\tau''$ stretch across $B'_1$, and $B_1'\cap(f\tau'\cup f\tau'')=B_1'\cap B^r(\overline 012{\bf 1})$.  Thus we see that $B_1^b$, as defined by Figure~5.5, is equal to $B_1'\cap f B_2'$.   Since $f\tau'$ and $f\tau''$ also stretch across $B'_0$, we see that $B_0^b=B_0'\cap fB_2'$.  Mapping this forward, we see that $B^a_j$, as defined by Figure~5.5, is equal to $B_j'\cap fB_0'$, for $j=0,1$.  
\medskip

\epsfysize=1.3in
\centerline{ \epsfbox{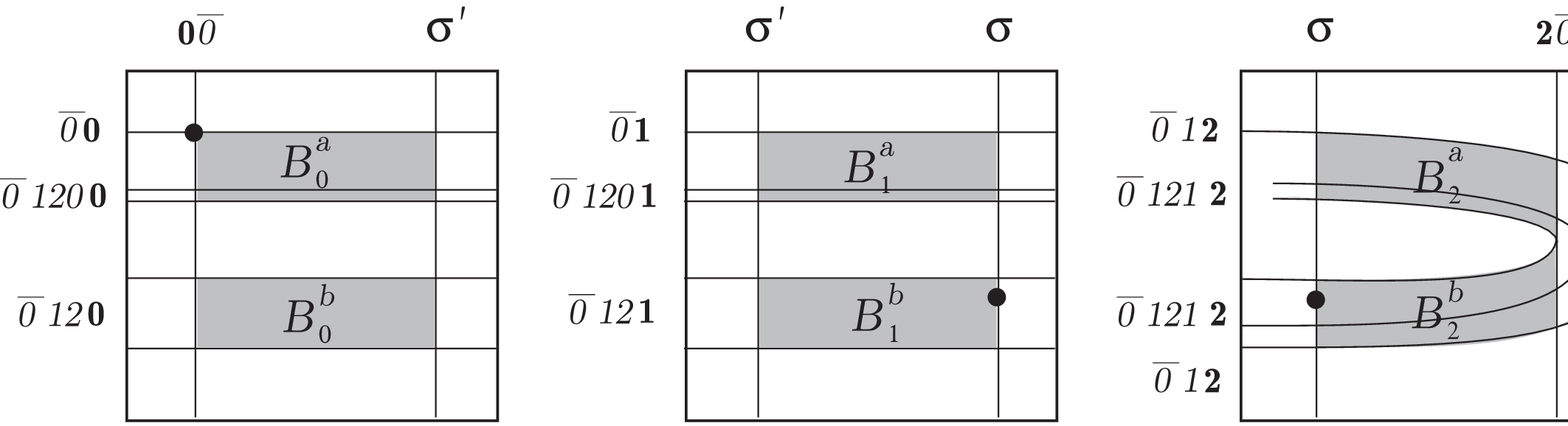}  }

\centerline{Figure 5.5.  Boxes $B_j^a$ and $B_j^b$, $0\le j\le 2$.  }

\medskip

We have seen that the vertical boundaries of $B_1^b$ are sub-arcs of $f\tau'$ and $f\tau''$ which have degree 1 in $B_1'$.  Mapping forward by $f_{1,2}$, we obtain two degree 2 arcs in $B_2'$.   By Figure~5.2, these sub-arcs are nested: there will be an inner one, which we call $\nu'$ and an outer one, which we call $\nu''$.  Let us denote their complexifications  by $\tilde\nu'$ and $\tilde\nu''$.  By Theorem 3.1, $\tilde\nu'$ and $\tilde\nu''$ are distinct complex disks of degree 2.    Counting multiplicity, we see that each disk $\tilde\nu'$ and $\tilde\nu''$ intersects $B({\bf 2}\overline 0)$ in 2 points.  Thus each of the curves $\nu'$ and $\nu''$ can intersect $B^r(\overline0{\bf 2})$ in at most 2 points.  By construction, we have that $S=(\nu'\cup\nu'')\cap B^r(\overline0{\bf 2})$.  Thus the inner curve $\nu'$ must intersect $B^r(\overline0{\bf 2})$ in 1 point if $\#S=3$, and 2 points if $\#S=4$.  In both cases, $\nu'$ contains  2 sub-arcs, which serve as the inner boundaries of $B_2^{a/b}$.

By construction, we see that $\bigcup_j(B_j^a\cup B_j^b)$ contains the union of the sets $B_0'\cap(B(0{\bf 0})\cup B(2{\bf 0}))$, $B_1'\cap(B(0{\bf 1})\cup B(2{\bf 1}))$, and $B_2'\cap B(1{\bf  2})$; and these sets contain $K_\infty$.  \qed

Let us end this section with a dichotomy on entropy.
\proclaim Theorem 5.6.  Let $f$ be a real H\'enon map with a real, 3-box system, and let $f^r:=f|_{{\bf R}^2}$ denote the restriction.  With $S$  as above, the real entropy $h(f^r)$ satisfies:  
\item{}  If $\#S\le 2$, then $h(f^r)<\log2$.
\item{}  If $\#S>2$, then $h(f^r)=\log2$. 

\noindent {\it Proof.}  If $h(f^r)=\log 2$, then by [BLS] all (complex) intersections between $W^s(p_0)$ and $W^u(p_0)$ occur inside ${\bf R}^2$.  The index $\overline 01212\overline0$ has degree 4, so the number of (complex) intersections between $B(\overline 0121{\bf 2})$ and $B({\bf 2}\overline0)$ is 4.  Thus we cannot have $\#S=0$, for otherwise there would be complex (non-real) intersections between $B(\overline0121{\bf 2})\subset W^u(p_0)$ and $B({\bf 2}\overline0)\subset W^s(p_0)$.

Now if $\#S\ge1$, then by Lemma 5.4, $(*)$ holds.  By the discussion before Lemma 5.2, we see that $B(\overline01{\bf 2})$ is a degree 2 curve opening to the left and intersecting $B({\bf 2}\overline0)$.  Thus $B(\overline 12{\bf 1})$ consists of two arcs of degree 1, and so $B(\overline0121{\bf 2})$ consists of two nested arcs of degree 2. By the property of real disks, each of these nested arcs is contained in a complex disk, which must have degree 2.  Thus each of the complex disks must intersect $B({\bf 2}\overline0)$ with multiplicity 2, and by the maximal entropy property, the real arcs must also intersect  $B^r({\bf 2}\overline0)$.  The arcs are coming in from the left, and they are nested; the inner arc intersects $B^r({\bf 2}\overline0)$ in at least one point, so the outer arc must intersect in two points.  This shows that if the entropy is $\log 2$, then $\#S\ge3$.

For the converse, we suppose that $\#S\ge3$.  Then we consider the boxes $B_j^{a/b}$ constructed in Lemma 5.5 and let $B^a=\bigcup_j B_j^a$ and $B^b=\bigcup_j B_j^b$.  We see that $f$ maps the boxes $B^{a/b}$ across themselves in the way that corresponds to a full 2-shift. 
 Consider the coding of a point $p\in K_\infty$ on the symbols $a$ and $b$, given by the itinerary.  If there is a point of tangency we remove its orbit.  Thus the coding will be well-defined outside of a countable set, giving a mapping onto a full shift $\Sigma_{a,b}$ on two symbols, minus countably many points.  Thus $(f^r,K_\infty)$ has entropy at least $\log 2$.  Since this is the maximum entropy possible for a quadratic H\'enon map, we see that $h(f^r)=\log2$.  \qed

\noindent {\bf 6.  Expansion. }  In the previous section, we found it useful to trim down the real boxes in a dynamically meaningful way.  In this section we will trim down the complex boxes.  This will allow us to control degrees of unstable disks, which in turn will allow us to show uniform expansion for our map.  Let us write $W^{s/u}_r:=W^{s/u}(p_0)\cap {\bf R}^2$. The following result describes what we can achieve with real boxes.
\proclaim Lemma 6.1.    If $\#S=4$, then for $0\le j\le 2$, each connected component of $W^u_r\cap B_j'$ is a horizontal arc of degree 1.  If $\#S=3$, then:
\item{$(i)$} The connected components of $W^u_r\cap B_1'$ are horizontal arcs of degree 1.
\item{$(ii)$}  All components of $W^u_r\cap B_2'$ have degree 1 except for one arc, written $\nu$, which has degree 2 and is tangent to $B^r({\bf 2}\overline0)$.
\item{$(iii)$}  All components of $W^u_r\cap B_0'$ have degree 1 except for the forward images under $f_{0,0}$ of the degree 2 arc $f_{2,0}\nu$.

\noindent{\it Proof. }  Suppose first that $\gamma$ is a horizontal arc of degree 1 in $B_0$ or $B_2$.  Then by the properties of the 3-box system, $f\gamma$ intersects $B_0$ and $B_2$ in horizontal arcs of degree 1.  The same holds with $B_i$ replaced by $B_i'$.  Next, we will consider the cases $\#S$ is 3 or 4.  For this we will make reference to Figures~5.4 and~5.5, which illustrate the results obtained in \S5.

Let $\gamma_0=B^r(\overline0{\bf 0})$ denote the local unstable manifold of $p_0$, so $W^u_r=\bigcup_{n\ge0}f^n\gamma_0$.   We will show by induction on $n$  that $f^n\gamma_0$ intersects $B_j'$ in real, horizontal arcs of degree 1, unless we are in the exceptional cases of $(ii)$ or $(iii)$.  By the previous paragraph, we need consider only horizontal arcs of degree 1 in $B_1$.   By Lemma 5.5, we see that $f\gamma$ intersects $B_2'$ in two arcs, one in $B_2^a$ and one in $B_2^b$.   If $\#S=4$, these arcs are disjoint, which completes the proof of that case, as well as the proof of $(i)$.  If $\#S=3$, then there is exactly one arc, the ``inner'' arc of the boundary of $B_2^a\cup B_2^b$, which does not break into two pieces.  This arc, which we will call $\nu$, is tangent to $B^r({\bf 2}\overline0)$, and verifies $(ii)$.  

We map $\nu$ forward to $B_0'$ and obtain an arc $f_{2,0}\nu$ of degree 2 in $B_0'$ which is tangent to the local stable manifold $B({\bf 0}\overline0)$.  Now we can continue to map $f_{2,0}\nu$ forward.  If we apply $f_{0,0}$ any number of times we obtain another arc which, because of the tangency, has degree 2 in $B_0'$.  If we apply $f_{0,1}$, then the arc is mapped to a pair of degree 1 arcs, and we are back in the case in the first paragraph above.  \qed

Now let us trim down the complex boxes.  We use the notation $0^n$ for the sequence $0\cdots0$.  The sets $B({\bf 0}0^n)$ are neighborhoods which shrink to the local stable manifold $B({\bf 0}\overline0)$ as $n\to\infty$.  Given $n$, we define
$$B_0^+:=B({\bf 0}0^n)  \cup  \bigcup_{0\le j\le n-1} B({\bf 0}0^j 1), \ \  B_1^+:=B_1\ \ {\rm and}\ \ B_2^+:=f^{-1}B_0^+.$$
Thus $B_j^+\subset B_j$, and we may shrink $B_j^+$ slightly so that the mappings strictly ``overflow,'' and we have a system $(\cB^+,\cG)$ of crossed mappings.
\proclaim Lemma 6.2.  If $\#S=4$, then we may choose $n$ sufficiently large that for each $j=0$, $1$ or $2$, the connected components of $W^u(p_0)\cap B_j^+$ are degree 1 horizontal disks.

\noindent{\it Proof. }  We saw in Lemma 6.1 that the components of $W^u_r\cap B_0'$ are degree 1 arcs.  Now $B({\bf 0}0^n)\cap {\bf R}^2$ shrinks down to the left hand boundary of $B_0'$, and the sets $B({\bf 0}0^j1)\cap {\bf R}^2$ are to the right of it, we may choose $n$ sufficiently large that the components of $W^u_r\cap B_0^+$ are also degree 1 arcs.  By $\S3$, then the components of $W^u(p_0)\cap B_0^+$ are horizontal disks of degree 1.  We get the statement for $W^u(p_0)\cap B_2^+$ by mapping backwards under $f$.  Finally the case $W^u(p_0)\cap B_1^+$ follows because $f_{0,1}$ and $f_{2,1}$ both have degree 1.   \qed

Now suppose that $\#S=4$.  Since the components of $W^u(p_0)\cap B_j^+$ all have degree 1, we may write them as graphs of (bounded) analytic functions.  Now let us define $\cW^u_j$ to be the set of all (degree one) disks in $B_j^+$ obtained as normal limits of these functions.  Thus the disks of $\cW^u_j$ contain all points of $K_\infty\cap B_j'$.  For $p\in K_\infty\cap B_j'$ we let $D_p$ denote the disk containing $p$.  Let $E^u_p\subset{\bf C}^2$ be the tangent space to  $D_p$, and for $v\in E^u_p$, let $s_p^u(v)$ be the Poincar\'e metric of $D_p$ applied to the vector $v$ at $p$.  

By the definition of a system of crossed mappings,  $f(D_p)$ overflows the disk $D_{fp}$ by a uniform amount.  Thus we have expansion with respect to the Poincar\'e metric:

\proclaim Lemma 6.3.  If $\#S=4$, then the metric $s_p^u$ is uniformly expanded under $f$.

By a similar argument we have expansion under $f^{-1}$.   We start by constructing boxes
$$B_0^-:=B(0^n{\bf 0})  \cup  \bigcup_{0\le j\le n-1} B(10^j{\bf 0} ), \ \   B_1^-:=fB_1 \ \ {\rm and}\ \ B_2^-:=B_2.$$
The only difference is that we need to see how vertical arcs in $B_j'$ map under $f^{-1}$.  The only map which is not of degree 1 is $(f^{-1})_{2,1}$, so we need only check what happens to vertical arcs in $B_2'$.   Figure~6.1 shows the box $B_2'$, which is bounded above and below by curves $\tau'$ and $\tau''$ (cf.\ Figure 5.4).  If $\#S\ge3$, then (cf.\ Figure~5.5), there are curves $\nu'$ and $\nu''$ which are arcs in $B(\overline0121{\bf 2})$ and which are contained in the horizontal boundary of $\partial(B_2^a\cup B_2^b)$.  If $\#S=4$, then $\nu'$ and $\nu''$ are disjoint.

\epsfysize=1.3in
\centerline{ \epsfbox{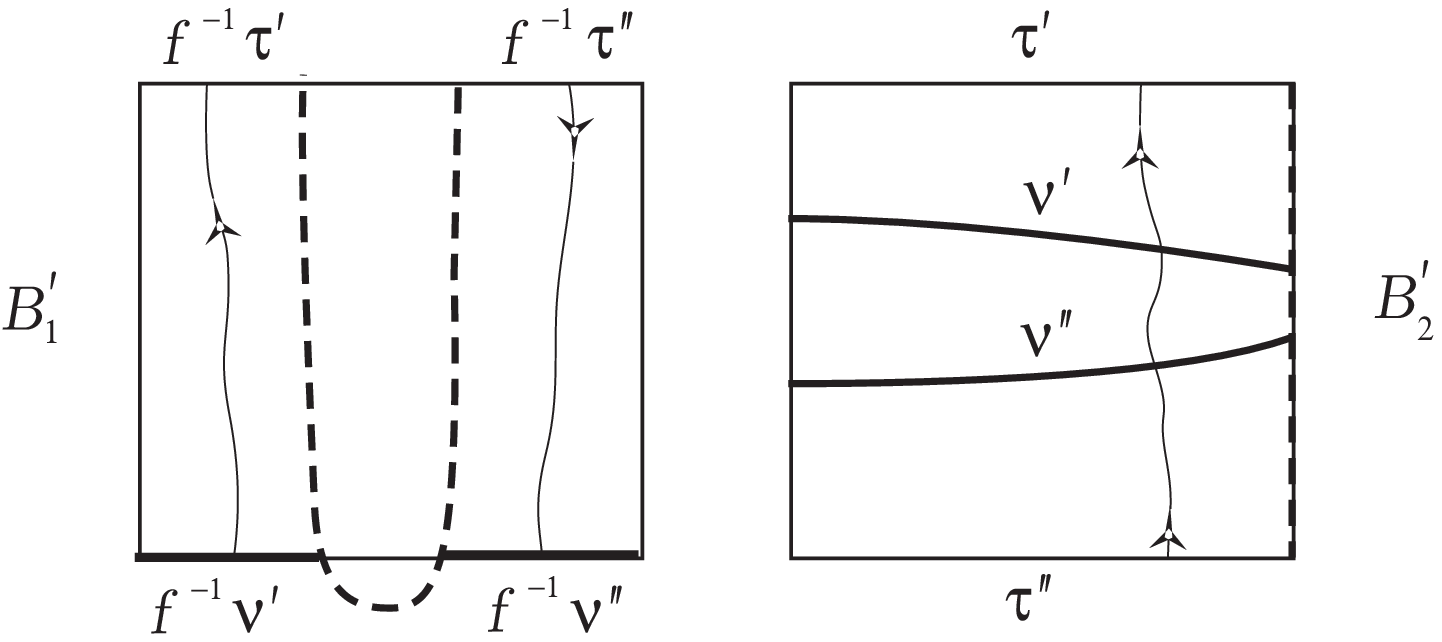}  }

\centerline{Figure 6.1.  Vertical arcs split under $f^{-1}:B_2'\to B_1$.  }

\medskip

The horizontal portions of $B_1'$ are contained in $W^u(p_0)$, as are the arcs $\tau'$, $\tau''$, $\nu'$ and $\nu''$.  Figure~6.1 shows that under $f^{-1}$ the arcs $\tau'$ and $\tau''$ are mapped to the upper portion of the boundary, and $\nu'$ and $\nu''$ are mapped to the lower portion of $\partial B_1'$.  The right hand boundary of $B_2'$ (which is dashed in Figure~6.1) is contained in $W^s(p_0)$, and the condition that $\#S=4$ means that $f^{-1}$ maps it to $B_1'$ as shown.  Thus we see that a vertical arc in $B_2'$ will be mapped to a curve that splits in $B_1'$.  We may apply Lemma 6.3 to $f^{-1}$, then to obtain the following:

\proclaim Theorem 6.4.  If $\#S=4$, then $f$ is hyperbolic on $K_\infty$.  Further, $(f,K_\infty)$ is conjugate to the full 2-shift.

\noindent{\it Proof.}  If $\#S=4$, we may split our 3-box system into the system shown in Figure 5.5.  In this system, the horizontal and vertical disks all have degree 1, so our map is hyperbolic on $K_\infty$. 
 Further, the sets $B^a:=B_0^a\cup B_1^a\cup B_2^a$ and $B^b:=B_0^b\cup B_1^b\cup B_2^b$ give a Markov partition.  The itinerary coding with respect to the partition $\{B^a, B^b\}$ gives a conjugacy to the 2-shift.  \qed

In case $\#S=3$, we may repeat many of the same arguments.  In this case, we get the following:
\proclaim Theorem 6.5.  If $\#S=3$, then we may choose $n$ sufficiently large  in the construction of $B^+_j$ so that for $j=0$, $1$ or $2$ each component of $B^+_j\cap W^u(p_0)$ has degree at most 2.  Similarly, each component of $B^-_j\cap W^s(p_0)$ has degree at most 2.

\bigskip\noindent{\bf 7.  External rays. }  
We continue with the hypothesis that $\#S\ge3$, as well as the notation from \S6.  This section will end with a proof of Theorem 2. 


  \medskip
\epsfysize=2.1in
\centerline{ \epsfbox{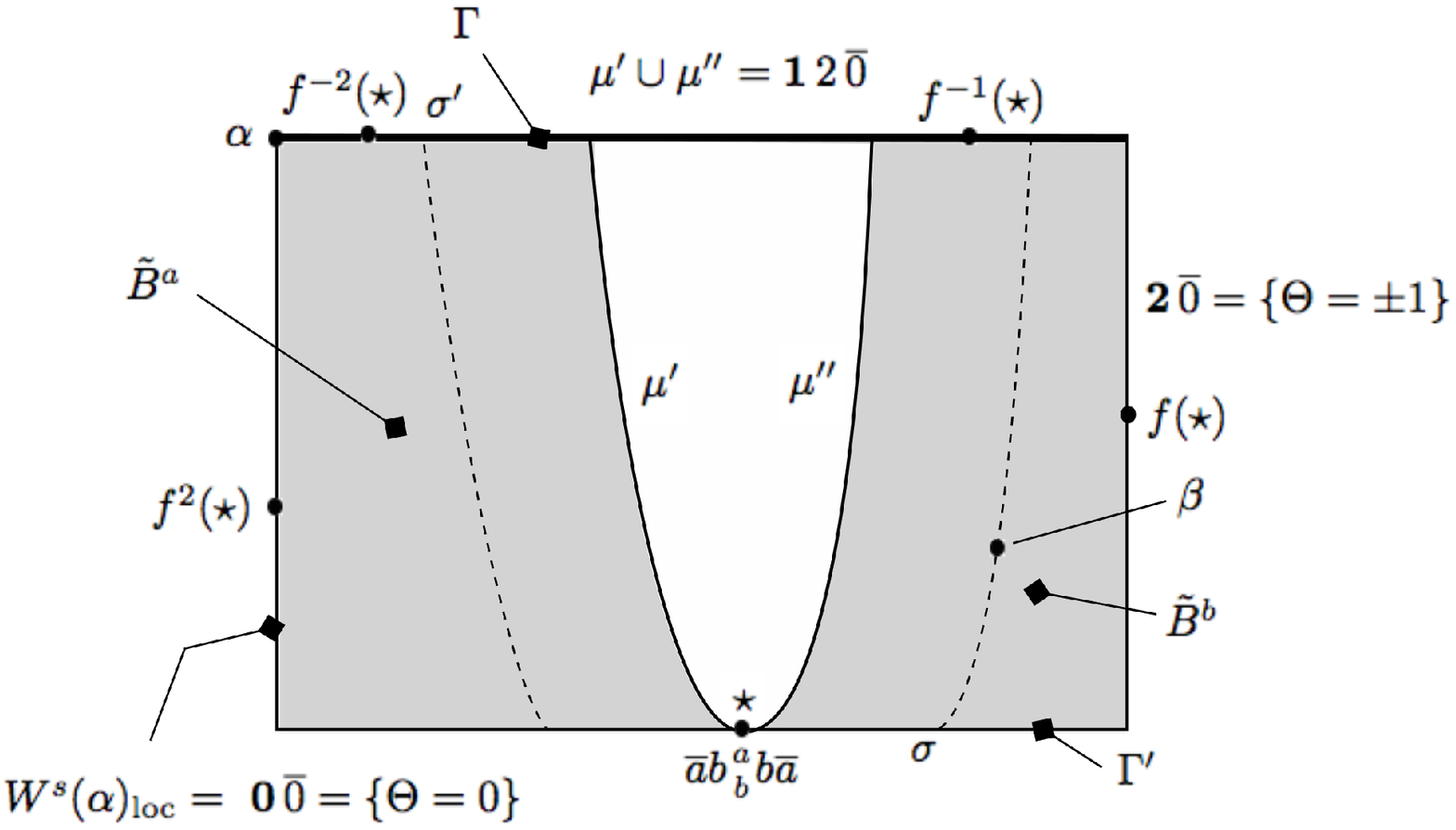}  }

\centerline{Figure 7.1.  Scheme for Coding.  }

\medskip

Let us consider a new box $\tilde B\subset{\bf R}^2$, shown in Figure 7.1, which is closely related to the union of boxes in Figure~5.5. The figure shows a ``rectangle'' whose boundaries are arcs of stable and unstable manifolds.    The fixed points are $\alpha$ and $\beta$.  The top and bottom portions of $\partial\tilde B$, labeled $\Gamma$ and $\Gamma'$, are arcs inside $W^u(\alpha)$.  The right and left portions of $\partial\tilde B$ are the arcs ${\bf 0}\overline 0=W^s_{\rm loc}(\alpha)$ and ${\bf 2}\overline0\subset W^s(\alpha)$.  The arcs $\sigma$ and $\sigma'$ are shown to make it easier to see the connection with Figure~5.5.

In Figure 5.5, we have a pair of sets $B^a=B^a_0\cup B^a_1\cup B^a_2$  and $B^b=B^b_0\cup B^b_1\cup B^b_2$.   This is almost a partition of $J$, in the sense that $B^a\cap B^b$ is a single point, and $J\subset B^a\cup B^b$.  
  In Figure~7.1 we have indicated the curve $\mu'\cup\mu''={\bf 1}2\overline0=f^{-1}_{12}({\bf 2}\overline0)$.  By the condition that $\#S=3$, there is a point of intersection $\star:=\mu'\cap\mu''$.  The image under $f^{-1}$ of the pair $\{B^a, B^b\}$ is the pair $\{\tilde B^a, \tilde B^b\}$.

 Let us write $J':=J-\bigcup_{n\in{\bf Z}}f^n(\star)$.  We define an itinerary coding map $s:J'\to \{a,b\}^{\bf Z}$
  by setting $s(p)=(s_j)_{j\in{\bf Z}}$, where $s_j=a$ if $f^j(p)\in B_a$ and $s_j=b$ otherwise.  Note that since $f$ maps $\tilde B^{a/b}$  to $B^{a/b}$, the coding with respect to one of these partitions differs only by a  shift from the coding with respect to the other partition.  Thus we will drop the tilde and just write $B^{a/b}$ for the partition in Figure~7.1.  The coding map is well defined on $J'$, but if we observe the part of the orbit of $\star$ shown in Figure 7.1, we see that the point $\star$ has two possible codings: $\overline a bab\overline a$ and $\overline a bbb\overline a$.  In order to deal with this we will pass to the quotient $\{a,b\}^{\bf Z}/\sim$, where we identify the two points $\overline a bab\overline a \sim \overline a bbb\overline a$, as well as all the pairs obtained by applying the shift.  This defines a relation which is closed when viewed as a subset of the product space. The quotient $\{a,b\}^{\bf Z}/\sim$ becomes a compact Hausdorff space when given the quotient topology (see [W] p.\ 128).  
  
  \proclaim Proposition 7.1.    $s$ extends continuously to a semi-conjugacy $s:J\to \{a,b\}^{\bf Z}/\sim$.  
  
  \noindent{\it Proof.} The continuity of $s$ at points that are not on the orbit of $\star$ is straightforward thus  it suffices to show that $s$ is continuous at $\star$.  A neighborhood basis of $s(\star)$ in the quotient topology is given by the cylinder sets $\{*a^Nb?ba^N*\}$
     where `$?$' denotes $a$ or $b$, and $*$ denotes arbitrary half-infinite sequences. 
     With this notation, we have $s(\star)=\overline ab? b\overline a$.   The orbit of $\star$ approaches $\alpha$ in both forward and backward time as is shown in Figure 7.1.  Given $N$, we may choose a neighborhood $U$ of $\star$ such that $f^jp\in B^a$ for all $p\in U\cap J$ for all $-N\le j\le N$ with $j\ne-1,0,1$.  Thus $s(p)$ will lie in the cylinder set about $s(\star)$.  It follows that $s$ is continuous at $\star$.    \qed
  
Since $s$ is a continuous mapping to a compact Hausdorff space, it remains only to  show that $s$ is a bijection, and it will then follow that $s$ is a homeomorphism and a conjugacy.  
  
Now we define external rays as in [BS7].  Let $G^+$ denote the Green function.  For a  disk $D:=D^u$ inside an unstable manifold, $G:=G^+|_{D}$ will be harmonic and strictly positive on the set $D-K^+$.  The {\it external rays} are exactly the gradient lines of $G$ on the sets $D-K^+$.
Since $J$ is real, we may take $D$ to be invariant under complex conjugation.  
Let $\gamma$ be an external ray inside $D-{\bf R}^2$.  If we follow $\gamma$ in the direction of decreasing $G^+$, then by Theorem 7.2 below, it will land  at a point of $D\cap {\bf R}^2$. 
We call this point $e(\gamma)\in D\cap {\bf R}^2$.     We have shown that $D\cap {\bf R}^2$ is an interval in Theorem 4.2, so $(D\cap{\bf R}^2)-J$ is a union of open intervals $I$.  We say that $q\in J$ is {\it exposed} if it is a boundary point of one of these open intervals.   Since the Green function $G^+$ of $J\cap D$ is continuous it follows that $J\cap D$ has no isolated points.  Thus $q$ can be the endpoint of only one interval, which we denote as $I_q$, and we may think of the point as {\it one-sided}.  We will use ${\cal N}$ to denote the set of gradient lines/external rays that land at non-exposed points, and ${\cal C}$ for the ones that land at critical points.

If $G$ is harmonic in a neighborhood of a point $\omega_0$ then, since $G$ is the real part of a holomorphic function, we may write
$$G(z) = G(\omega_0)+ \Re\left(c (z-\omega_0)^m\right) + \cdots\eqno{(*)}$$ 
for some nonzero constant $c$.  If $m\ge2$, $\omega_0$ is a {\it critical point} and we call $m$ the {\it order} of the critical point. If $m=2$ we say the critical point is {\it simple}.  If $D^u$ is an unstable disk, then we will consider $G:=G^+|_{D^u}$.  For $\lambda>0$ small, the function $\lambda^{-1}G$ and any connected component $\Omega$ of $\{G<\lambda\}$ will satisfy the hypotheses of Theorem 7.2.  The following result gives the landing of external rays:
\proclaim Theorem 7.2.  Let $\Omega\subset{\bf C}$ be a bounded open domain which is simply connected, smooth, bounded and  invariant under complex conjugation.  Let  $E\subset\Omega\cap {\bf R}$ be compact.  Let $G$ be a continuous function on $\bar \Omega$ which is equal to 1 on $\partial\Omega$, is equal to $0$ on $E$ and is harmonic on $\Omega-E$.  Then we have:
\item {($i$)}  $G$ has no critical points on $\Omega-{\bf R}$.
\item{($ii$)}  Each proper complementary interval contains exactly one critical point, which is simple.
\item {($iii$)}  For every $p\in\partial\Omega$ there is a point $e(p)\in {\bf R}$ and there is a gradient curve of $G$, $\gamma\subset\Omega$, such that $\gamma$ makes a continuous arc connecting $p$ to $e(p)$.
\item{($iv$)}  Each non-exposed point of $E$ and each critical point is the landing point of a unique conjugate pair of rays.
\item{($v$)}   
The landing map is continuous at all points of ${\cal N}$ and discontinuous at all points of ${\cal C}$.   $\partial \Omega-{\bf R}$ is equal to the set of initial points of  ${\cal N}\cup{\cal C}$.

\noindent {\it Proof. }  If $\omega_0\in\Omega-E$ is a critical point of $G$ of order $m$, then near $\omega_0$ there are $m$ gradient curves leaving $\omega_0$ along which $G$ is increasing, and these are separated by $m$ gradient curves along which $G$ is decreasing.  If $\omega_0$ is critical, then $m$ is at least 2.  Thus at every critical point $\omega_0$, there are at least two gradient lines $\gamma_0'$ and $\gamma_0''$ along which  $G$ is decreasing.  Now suppose that $\omega_0\in\Omega-{\bf R}$ and choose one of the decreasing gradient lines, say $\gamma_0'$.  Since $\{G\ge\lambda>0\}$ contains only finitely many critical points, we may follow $\gamma_0'$ until it ends at another critical point, or otherwise we may follow it and obtain a curve with $\inf_{\gamma_ 0}G=0$.  If $\gamma_0'$ ends in a critical point $\omega_1$, we may choose a gradient line $\gamma_1'$ exiting $\omega_1$ and continue in this way so that we have a curve $\gamma':=\bigcup_j\gamma'_j$ along which we have $\inf_{\gamma'}G=0$.   Similarly, there is a second curve $\gamma_0''$ emanating from $\omega$ which can be extended to $\gamma''$, with $\inf_{\gamma''}G=0$.  Now since $\omega_0\notin{\bf R}$, we have a conjugate critical point $\bar\omega_0$ and curves $\bar\gamma'$ and $\bar\gamma''$.  Now let $\Gamma$ denote the closure of $\gamma'\cup\gamma''\cup\bar\gamma'\cup\bar\gamma''$.  If follows that there is a bounded component $\Omega_0\subset\Omega$ with $\partial\Omega_0\subset\Gamma$.  Since the boundary of $\Omega_0$ consists of gradient lines, we conclude that $\max_{\bar\Omega_0}G=G(\omega_0)$.  This is a contradiction because there is also a gradient line along which $G$ increases, and this is between $\gamma_0'$ and $\gamma_0''$, so it enters $\Omega_0$.  This proves $(i)$.

Now let us suppose that $\omega_0\in {\bf R}\cap (\Omega-E)$ is a critical point.  By conjugation symmetry, we have $G(z)=G(\bar z)$.  Thus the constant $c$ in $(*)$ must be real.  If $m$ is even, then the real arcs emanating from $\omega$ on either side  are gradient lines of the same kind: either $G$ is increasing (or decreasing) in both directions and thus $G|_{\bf R}$ either has a strict local minimum or maximum at $\omega_0$.  If $m$ is odd, then $G|_{\bf R}$ strictly increases to one side of $\omega_0$ and strictly decreases on the other.

Let $I_{\omega_0}$ denote the interval of ${\bf R}-E$ which contains $\omega_0$.  We know that $G>0$ on this interval, and since it is proper, $G$ vanishes on both endpoints.  Thus there must be at least one maximum point, which must be critical and a strict local maximum.  Let us suppose that there is a second critical point, say $\omega_1$.  Suppose that the order $m$ of $\omega_1$ is odd.  In this case we may suppose that $G|_{\bf R}$ is strictly decreasing on the arc $[\omega_0,\omega_1]$.  Since the order is odd, it follows that $m\ge3$, and so there is a gradient arc $\tau$ emanating from $\omega_1$ into the direction of increasing $Im(z)$.  We may extend $\tau$ (resp.\ $\gamma$) to a gradient arc connecting $\omega_1$ (resp.\ $\omega_0$) to $\partial\Omega$.  Thus $G$ is a harmonic in a region $\Omega_0$ with $\partial\Omega_0\subset \gamma\cup\tau\cup I_{\omega_0} \cup\partial\Omega$.  By the Minimum Principle,  the minimum of $G|_{\bar \Omega_0}$ is attained at the point $\omega_1$.  However, since $\omega_1$ is a point of order $3$ or more, there is a gradient line $\nu$ which emanates from $\omega_1$, entering $\Omega_0$ so that  $G$ decreases along this arc.  This contradiction shows that $\omega_1$ of odd order cannot exist.  

A similar argument involving gradient arcs of increasing $G$, emanating from $\omega_0$ in the direction of increasing $Im(z)$, shows that $\omega_0$ must be a simple critical point.  If it is not simple, then it must be of even order at least 4.  The two rays leaving $\omega_0$ inside ${\bf R}$ are decreasing, and so there is at least one ray $\eta$  which is decreasing for $G$ and leaving $\omega_0$ in a direction that enters the upper half plane.  It follows that there are gradient rays $\gamma'$ and $\gamma''$ on either side of $\eta$ which enter the upper half plane.  $G$ increases along these rays, and we may follow them until we reach $\partial\Omega$.  Let $\Omega'$ denote the domain bounded by $\gamma'\cup\gamma''\cup\partial\Omega$.  It follows that $G(\omega_0)$ is the minimum value of $G$ on $\partial\Omega'$.  On the other hand, $\eta$ enters $\Omega'$, and $G$ decreases along $\eta$, which contradicts the Miminum Principle for harmonic functions.   This establishes $(ii)$.

To prove $(iii)$, we let $\gamma$ be a gradient curve emanating from $p\in\partial\Omega$.  There are no critical points in $\Omega-{\bf R}$, so it has unlimited continuation in that region.  We consider the set ${\bf cl}(\gamma)$, which is the set of cluster points of the curve.  This must be a closed, connected set and is contained in ${\bf R}$.  It suffices to show that it consists of a single point of $E\cup \{{\rm Critical\  points}\}$.  Otherwise, suppose that ${\bf cl}(\gamma)$ is an interval $I\subset{\bf R}$.  First, suppose that $I$ contains a critical point $q$.  Then $G$ has the form $(*)$ near $q$, so $\gamma$ must coincide with the gradient arc emanating from $q$.  Thus $\gamma$ lands at $q$, so $e(p)=q$.  Otherwise, since $\gamma$ is a gradient line of decreasing $G$, and $G|_E=0$, we see that we must have $I\subset E$.  But consider a point $r$ in the interior of $I$.  If we extend $G$ over $I$ to the lower half plane by setting it equal to  $-G(\bar z)$, then by the Reflection Principle, this extension is harmonic in a neighborhood of $I$.  Thus we have $G=\alpha(z)|y|$ for some smooth, non-vanishing function $\alpha(z)$, and each gradient line of such a function lands at a single point.  Since $\gamma$ enters the set where this behavior occurs the cluster set of $\gamma$ is a single point $e(p)$.  

To prove $(iv)$, we start by considering $p_1, p_2\in\Omega\cap \{Im(z)>0\}$.  First we show that $e(p_1)\ne e(p_2)$.   We know that the gradient lines $\gamma_{p_1}$ and $\gamma_{p_2}$ are disjoint, so  if $e(p_1)=e(p_2)$, then the arc $\overline{p_1,p_2}$ inside $\partial\Omega$, together with $\gamma_{p_1}\cup\gamma_{p_2}$ bound a simply connected domain $D$.  Let $H$ denote the harmonic conjugate of $G$, so that $G+iH$ gives a conformal map of $D$ into the strip $\{0<Re(z)<1\}$.  Since $\gamma_{p_1}$ and $\gamma_{p_2}$ are gradient lines of $G$, it follows that they are level sets for $H$.  Thus $D$ is mapped conformally onto the rectangle $\{0<Re(z)<1, a_1<Im(z)<a_2\}$.  This is not possible, for otherwise the inverse $(G+iH)^{-1}$ would be constant on the interval $[a_1i,a_2i]$ of the boundary of the rectangle.

Thus for each $p\in \partial\Omega$, we have that $p$ and $\bar p$ are the only points of $\partial\Omega$ which are mapped to the point $e(p)=e(\bar p)$.  As noted before, the map $p\mapsto e(p)$ is monotone.  Further, we claim that it is onto.  We have seen that for each critical point $c\in{\bf R}-E$ there is a gradient line (and thus a conjugate pair of gradient lines) landing at $c$.  By the uniqueness statement of the previous paragraph, we see that if $x\in I_c$ is an exposed point, then there can be no $p\in\partial\Omega$ which lands at $x$.  Now it follows by monotonicity that all the non-exposed points are landing points of some $p\in\partial\Omega$.    \qed

\noindent{\bf Remark.}  There is a connection between  Theorem 7.2, which concerns more general harmonic functions, and the dynamical properties of $G^+$ and ${\cal W}^u$ we have developed in \S6.  For instance, by Theorem 6.5, the connected components of $B_j^+\cap W^u(p_0)$ have degree $\le2$.  Critical points represent tangencies between $W^u(p_0)$ and the level sets of $G^+$.  This could be used to give an alternative proof of the statement in $(ii)$ that critical points are simple.

In the case $\#S=4$, Figure 7.1 will be different:  $\mu'$ and $\mu''$ do not intersect,  but they bound an interval $[\mu'\cap\Gamma',\mu''\cap \Gamma']$ inside $\Gamma'$.  This interval contains a critical point, which we may regard as $\star$.  The condition $\#S=3$ means that this interval has shrunk to a point, and the critical point has been lost and is now a point of $J$.

\medskip\noindent{\bf External angles.}   Now we will assign angles to external rays.  
 There is a unique holomorphic function $\varphi^+:\{ (x,y)\in{\bf C}^2:|x|>\max(R,|y|)\} \to {\bf C}$ with the property that $-\varphi^+(x,y)\sim x$ near infinity and $\varphi^+\circ f = (\varphi^+)^2$.   Further, $G^+ = \log|\varphi^+|$, so $\varphi^+$ has a  holomorphic continuation along any curve in the region  ${\bf C}^2-K^+$, since $G^+$ is harmonic there. (The continuation may depend on the choice of the curve).  Since $Arg(\varphi^+)$ and $G^+$ are harmonic conjugates, the gradient curves of $G^+$ are level sets for the argument of $\varphi^+$.  If we continue $\varphi^+$ in from the points of ${\bf R}^2$ with $\Re(x)\ll0$, then we may continue it along the real locus until we reach the left hand boundary $W^s(\alpha)_{\rm loc}$ of $\tilde B$, and $\varphi^+$ takes the value $+1$ there.  Thus we choose the argument of $\varphi^+$ to be 0 along this curve.   In \S6 we showed that each $p\in J\cap B_j$ is contained in an unstable disk $D:=D^u_p$ which is horizontal in $B_j$. 
 Since $D$ is invariant under complex conjugation, we may use Theorem 4.1 and suppose that that $D-{\bf R}^2$ consists of two simply connected pieces. 
 Further,  since $D\cap K^+\subset{\bf R}^2$, we may suppose that the branch of $\varphi^+$ extends analytically  to $D-{\bf R}^2$ so that it approaches $+1$ on the left hand boundary of $\tilde B$.   It will then approach $-1$ at the right hand vertical boundary of $\tilde B$.

We will use the Principal Argument function, which takes values in the interval $[-\pi,\pi)$. We index the gradient lines in $D$ as: $\gamma_{\theta} :=\{Arg(\varphi^+)=\pi \theta \}$. 
(Note that our convention differs from the usual convention in one variable by a factor of 2.) Thus $D-{\bf R}^2$ is filled by the curves $\gamma_\theta$ for $-1<\theta<1$.  
  By Theorem 7.2 each $\gamma_\theta$
  will land either on a non-exposed point of $J\cap D^u_p$ or at a critical point of $G^+$ in $D^u_p\cap {\bf R}^2$.  If $p$ is not exposed, then there is a  conjugate pair of external rays  $\gamma_{\pm \theta}\subset D^u_p$, $0<\theta<1$, which land at $p$; and we set $\Theta(p)= \theta$, with $0< \theta< 1$.  Otherwise,  $\gamma_\theta$ lands at a critical point $c$ with $0<\theta<1$.  In this case, we let $I_c=(p_-,p_+)$ denote the associated complementary interval.  We define $\Theta(p_+)=\Theta(p_-)=\theta$.  
 It follows that we have $\Theta(f(p))=2\Theta(p)$ modulo~1.  By symmetry, we have $\Theta=1$ on the right hand boundary of $\tilde B$.  Thus we have the following:
\proclaim Corollary 7.3.  For each $0\le \theta\le 1$ and each $p\in J$ the set $D_p\cap\{\Theta=\theta\}$ is nonempty and consists of 1 or 2 points.

The level sets  $\{p\in J:\Theta(p)=\theta\}$ give a partition of $J$, and we describe these sets in more detail.  Let us work with the complex vertical (local stable) disk $V_0=B({\bf 0}\overline 0)$.  This is a vertical disk of degree 1 inside the complex box  $B_0$, so $V_0\cap J\subset\{\Theta=0\}$.  We iterate $V_0$ backwards, so that$f^{-n}V_0\cap J$ is the union of the sets $\{\Theta={j\over 2^n}\}$.   Let $V$ be a component of $f^{-n}V_0\cap B$,  so the argument of $\varphi^+$ is constant as we approach $V\cap J$ so it is contained in one of the sets $\{\Theta={j\over 2^n}\}$ for some $j$.   By Theorem 6.5,  $V$ is a vertical disk of degree $\le 2$. 
 Since the degrees of these disks are bounded by 2, we may take the limits of them and obtain vertical disks of degree $\le 2$.  Since the dyadic numbers are dense in $[0,1]$, we see that we may take limits so that for each $\theta$, there will be a vertical disk $V_\theta$ such that $\Theta=\theta$ on $V_\theta\cap J$.  Let us define $\tau_\theta:=V_\theta\cap \tilde B\subset{\bf R}^2$.  By \S5, we know that $\tau_\theta$ will have degree 2 only if $\tau_\theta$ intersects $\Gamma\cup\Gamma'$ in a point of the form $f^m(\star)$ for some $m\in{\bf Z}$. 
To summarize, we have:
\proclaim Lemma 7.4.  For each  $0\le \theta\le 1$ there is an  arc $\tau_\theta$ (or pair of arcs $\tau_{\theta^-}$ and $\tau_{\theta^+})$ which is vertical in $\tilde B$ with degree 1 and which has the property:  $\{p\in J:\Theta(p)=\theta\}\subset\tau_\theta$ (or $\{p\in J:\Theta(p)=\theta\}\subset\tau_{\theta^-}\cup\tau_{\theta^+}$). In other words,  $J\cap\tilde B$ is contained in the union of smooth (real) vertical curves  $\tau_\theta$, where $\Theta(p)=\theta$ for all $p\in J\cap \tau_\theta$.  Each $\tau_\theta$ is an arc of a stable manifold in ${\cal W}^u\cap{\bf R}^2$, and $\tau_\theta\cap J$ is nonempty.

The function $\Theta$
is closely related to the $s$-coding defined at the beginning of this section.  Namely for $s=(s_j)_{j\in{\bf Z}}$,  we define  $\Theta(s)=(\theta_j)_{j\in{\bf Z}}$, such that $\theta_j$ is the number whose binary coding is $.d_jd_{j+1}d_{j+2}\cdots$, where $d_j=1$ if $s_j=b$ and 0 otherwise.   

Let ${\cal S}'=\{a,b\}^{\bf Z} - \bigcup_{j\in{\bf Z}}f^j\{\overline a bab\overline a, \overline a bbb \overline a\}$.  
\proclaim Lemma 7.5.  For every $(s_j)_{j\in {\bf Z}}\in{\cal S}'$, there is a point $p\in J'$ such that $s(p)=(s_j)_{j\in {\bf Z}}$.  

\noindent{\it Proof.}  Let $s\in{\cal S}'$ be given, and let $(\theta_j)_{j\in{\bf Z}}:=\Theta(s)$.
For each $n$, choose a point $p_n\in \tau_{\theta_{-n}}$.   It follows that the $\Theta$-code for $p_n$ is the same as $(s_j)$ for all $j\ge -n$.   Since $J$ is compact, we may choose $p\in\tau_{\theta_{0}}$ to be a limit point of the sequence $f^n(p_n)\in\tau_{\theta_0}$.  If $p$ is in the $f$-orbit of $\star$, then it follows that $s\notin{\cal S}'$.  Otherwise, we have $s(p)=s$.  \qed

\noindent{\it Proof of Theorem 2.} 
By Lemma 7.5, we know that $s$ is surjective.  It remains to show that it is injective.  We have shown that $J'$ is stratified by a family of stable arcs $\tau_\theta$; each $p\in J'$ lies in a unique vertical arc $\tau_\theta$.  Now we wish to obtain a similar stratification in terms of unstable arcs.  For this, we recall that there is also a holomorphic function $\varphi^-:  \{ (x,y)\in{\bf C}^2:|y|>\max(R,|x|)\} \to {\bf C}$ with the property that $\varphi^-\sim y/b$ near infinity, and $\varphi^-\circ f^{-1} = (\varphi^-)^2$.  
For each $p\in J'$ and each stable disk $D^s$ containing $p$, we may consider the rays $\eta_{\vartheta}:=\{Arg(\varphi^-|_{D^s})=\pi\vartheta\}$, which are defined as the sets of $Arg(\varphi^-|_{D^s})$.

Now we consider the partition $\{B^a, B^b\}$ corresponding to Figure 5.5.  If $c=(c_j)_{j\in{\bf Z}}\in\{a,b\}^{\bf Z}/\sim$, then we may pass to a coding $(\vartheta_j)$, and we have the corresponding curve $\eta_\vartheta$.  Since we are intersecting horizontal and vertical curves of degree 1 in $\tilde B$, it follows that there is a unique point $p\in \tau_\theta\cap \eta_\vartheta$ with $s(p)=c$.  \qed

\medskip\centerline{\bf References}
\medskip
\item{[BLS]}  E. Bedford, M. Lyubich, and J. Smillie, Polynomial
diffeomorphisms of $\C^2$.  IV: The measure of maximal entropy and laminar
currents,  Invent.\ Math.\ 112, 77--125 (1993).

\item{[BS7]}  E. Bedford and J. Smillie,  Polynomial diffeomorphisms of
$\C^2$. VII: Hyperbolicity and external rays.  Ann.\ Sci.\ \'Ecole Norm.\ Sup.\ (4) 32 (1999), no.\ 4, 455--497.

\item{[BS8]}  E. Bedford and J. Smillie,  Polynomial diffeomorphisms of
$\C^2$.  VIII: Quasi-expansion.  American J. of Math., 124, 221--271, (2002).

\item{[BSh]} E. Bedford and J. Smillie, The H\'enon family: The complex horsehoe locus and real parameter space.  {\sl Complex dynamics}, 21--36, Contemp.\ Math., 396, Amer.\ Math.\ Soc., Providence, RI (2006).

\item{[BS$i$]} E. Bedford and J. Smillie, Real polynomial diffeomorphisms with
maximal entropy: Tangencies.  Annals of Math. 160 (2004), 1--25.

\item{[BS$ii$]} E. Bedford and J. Smillie, Real polynomial diffeomorphisms with maximal entropy. II. Small Jacobian. {\sl Ergodic Theory Dynam. Systems} 26 (2006), no. 5, 1259--1283. 

\item{[CLR1]} Y. Cao, S. Luzzatto and I. Rios,  Some non-hyperbolic systems with strictly non-zero Lyapunov exponents for all invariant measures:  horseshoes with internal tangencies, Discrete Cont.\ Dyn.\ Syst., 15 (2006), 61--71.

\item{[CLR2]} Y. Cao, S. Luzzatto and I. Rios,  The boundary of hyperbolicity for some H\'enon-like families, The boundary of hyperbolicity for H\'enon-like families. Ergodic Theory Dynam.\ Systems 28 (2008), no. 4, 1049--1080.  arxiv/math 0502235

\item{[DCH]}  A. De Carvalho and T. Hall, How to prune a horseshoe, Nonlinearity, 15 (2002), R19--R68.


\item{[DDS]}  T-C Dinh, R. Dujardin, and N. Sibony,  On the dynamics near infinity of some polynomial mappings in ${\bf C}^2$, Math.\ Ann.\ 333 (4) (2005), 703--739.

\item{[Du]}  R. Dujardin, H\'enon-like mappings in ${\bf C}^2$, Amer.\ J. Math.\ 126 (2) (2004), 439--472.

\item{[HS]} R. Hagiwara and  A. Shudo,  An algorithm to prune the area-preserving H\'enon map. {\sl J. Phys. A} 37 (2004), no. 44, 10521--10543. 

\item{[H1]}  U. Hoensch,  Horseshoe-type diffeomorphisms with a homoclinic tangency at the boundary of hyperbolicity. Thesis (Ph.D.)  Michigan State University. 2003.  57pp.

\item{[H2]} U. Hoensch, Some hyperbolicity results for H\'enon-like diffeomorphisms. 
Nonlinearity 21 (2008), no. 3, 587--611.

\item{[HO]}  J.H. Hubbard and R. Oberste-Vorth, H\'enon mappings in the complex domain.  II:~Projective and inductive limits of polynomials.  {\sl Real and Complex Dynamical Systems} (Hiller\o d, 1993), 89--132, NATO Adv.\ Sci.\ Inst.\ Ser.\ C Math.\ Phys.\ Sci., 464, Kluwer Acad. Publ., Dordrecht (1995).  



\item{[IS]} Y. Ishii and J. Smillie, Homotopy shadowing.   {\sl Amer. J. Math.} 132 (2010), no. 4, 987--1029. 


\item{[L]} C. Lipa,  Monodromy and H\'enon mappings, Thesis, Cornell University, 2009.

\item{[M]} J. Milnor,  Periodic orbits, external rays and the Mandelbrot set: an expository account. 
G\'eom\'etrie complexe et syst\`emes dynamiques (Orsay, 1995).
Ast\'erisque No. 261 (2000), xiii, 277--333. 

\item{[T]} H. Takahasi, Prevalence of non-uniform hyperbolicity at the first bifurcation of H\'enon-like families. 	arXiv:1308.4199 [math.DS].

\item{[W]} S. Willard, General Topology, Addison-Wesley, 1970.

\bigskip
\rightline{Eric Bedford}

\rightline{Indiana University}

\rightline{Bloomington, IN 47405}

\rightline{\tt bedford@indiana.edu}

\rightline{current address: Stony Brook University}

\bigskip
\rightline{John Smillie}

\rightline{Cornell University}

\rightline{Ithaca, NY 14853}

\rightline{\tt smillie@math.cornell.edu}

\rightline{current address: Warwick Math Institute}

\end
 
 \proclaim Theorem. If (*) holds, and if $\epsilon'>0$ is given, then there exist $\epsilon,\delta>0$ with the following properties.  If $z$ is an $\epsilon$-orbit for $(f,K_{\infty})$, then there is an $\epsilon'$-orbit $(s,p_{s})$ for the shift on $\Lambda$ such that  $\phi(s,p_{s})=z$.  Now suppose that $z'$ and $z''$ are two such $\epsilon$-orbits, and suppose that the corresponding $\epsilon'$-orbits are $(s',p'_{s'})$ and $(s'',p''_{s''})$.  If $z'$ and $z''$ are $\delta$-close, then $s'=s''$.

We will say that the imbedding $\phi$ of the crossed mapping system $(\cB,\cG)$ is a {\it realization of} $(f,K_{\infty})$ if $(*)$ holds and if the map $\pi:\Lambda\to\Sigma_{\cG}$ has finite fibers.